\documentclass[onefignum,onetabnum]{siamonline250211}

\usepackage{lipsum}
\usepackage{amsfonts}
\usepackage{graphicx}
\usepackage{epstopdf}
\usepackage{algorithmic}
\usepackage{booktabs}

\newcommand{\md}{\mathrm d}

\newcommand{\E}{\mathbb{E}}
\newcommand{\PP}{\mathbb{P}}
\newcommand{\QQ}{\mathbb{Q}}
\newcommand{\RR}{\mathbb{R}}
\newcommand{\KL}{\mathrm{KL}}
\newcommand{\Law}{\mathrm{Law}}
\ifpdf
  \DeclareGraphicsExtensions{.eps,.pdf,.png,.jpg}
\else
  \DeclareGraphicsExtensions{.eps}
\fi

\usepackage{enumitem}
\setlist[enumerate]{leftmargin=.5in}
\setlist[itemize]{leftmargin=.5in}

\newsiamremark{remark}{Remark}
\newsiamremark{hypothesis}{Hypothesis}
\crefname{hypothesis}{Hypothesis}{Hypotheses}
\newsiamthm{claim}{Claim}
\newsiamremark{fact}{Fact}
\crefname{fact}{Fact}{Facts}

\newtheorem{thm}{Theorem}[section]
\newtheorem{prop}{Proposition}[section]
\newtheorem{cor}{Corollary}[section]
\newtheorem{asmp}{Assumption}

\theoremstyle{definition}

\headers{Path inference under partial observations}{}

\title{Pathwise Learning of Stochastic Dynamical Systems with Partial Observations
\thanks{\funding{NTY’s work is supported in part by National Science Foundation grant DMS-2038118.}}

}

\author{Nicole Tianjiao Yang\thanks{Department of Mathematics, University of Tennessee, Knoxville, TN 
  (\email{nicole.yang@utk.edu}).}
}

\usepackage{amsopn}

\ifpdf
\hypersetup{
  pdftitle={Bayesian Path Filtering},
  pdfauthor={N.T. Yang}
}
\fi

\begin{document}

\maketitle

\begin{abstract}
The reconstruction and inference of stochastic dynamical systems from data is a fundamental task in inverse problems and statistical learning. While surrogate modeling advances computational methods to approximate these dynamics, standard approaches typically require high-fidelity training data. In many practical settings, however, the data are indirectly observed through noisy and nonlinear measurements. The challenge lies not only in approximating the coefficients of the SDEs, but in simultaneously inferring the posterior updates given the observations. In this work, we present an amortized path generation method to address these challenges and solve nonlinear stochastic filtering from noisy observations. We first derive a variational inference formulation that solves filtering distribution for a given noisy observation path. This leads to a controlled SDE representation in which the feedback control is identified through the score structure of a pathwise Zakai equation. Motivated by this representation, we construct a conditional generative model that learns, in an amortized manner over observation paths, to transport a prior latent path measure toward the corresponding posterior path measure. We demonstrate the method on nonlinear stochastic systems with multimodal posterior structure, chaotic dynamics, and sparse observations, showing that the learned conditional path generator enables uncertainty quantification for both filtering marginals and trajectory-dependent functionals.
\end{abstract}


\section{Introduction}

Many scientific and engineering disciplines rely on the inference of the stochastic dynamical systems from noisy, partial observations. This is a central topic in the field of data assimilation, with many applications in geosciences, neurosciences, epidemiology, robotics, and financial engineering \cite{asch2016data, law2015data, reich2015probabilistic, brajard2020combining, moye2018data, engbert2021sequential, thrun2002probabilistic, jasra2011inference, yang2025finding}. By integrating real-world measurements with dynamical models, the goal is to improve perform accurate system inference as well as downstream tasks such as control, prediction, uncertainty quantification, and decision making.

Sequential Monte Carlo (SMC) methods \cite{doucet2000sequential}, commonly known as particle filters and particle smoothers, are among the most widely used tools for nonlinear state estimation. They approximate intractable posteriors by propagating ensembles of particles and repeated reweighting and resampling at observation times. While powerful, these methods faces several challenges: (i) a large number of particles is required to avoid degeneracy, which causes severe difficulty in high dimensions. 
(ii) the resampling step introduces discontinuities of the likelihood as a function of model parameters. In contrast, Ensemble Kalman Filters (EnKFs) \cite{burgers1998analysis, evensen2003ensemble, bergemann2012ensemble} offer computational efficiency but are inconsistent with Bayes theorem for non-Gaussian models. Their posterior distribution does not converge to the true filtering distribution as the number of particles $N \rightarrow \infty$.

Above all, the aforementioned methods typically require the governing equation or the transition kernel of true states and will struggle when the underlying dynamics are unknown or partially specified \cite{kantas2015particle}. Meanwhile, there has been increasing interests in data-driven methods to learn dynamical systems from time series data \cite{yu2024learning, karniadakis2021physics, park2024dynamical}, where the fidelity of training data is critical as well. We ask: Will the surrogate model robust under errors of the data? How can one extract the correct information from the noisy observations?

In this paper, we aim to develop an efficient data-driven method to perform simultaneous SDE inference and path estimation to answer these questions. Rather than updating the empirical marginal posterior density at each observation time in classical data assimilation approaches, we aim to learn a generative representation of the posterior path measure that is amortized with respect to the observation trajectory $y_{0:T}$. Concretely, we seek a map $y_{0:t}\ \mapsto\ \Pi_t^{\,y}$, where $\Pi_t^{\,y} \in \mathcal{P}(C([0,t];\mathbb{R}^n))$ denotes the conditional law of the signal path $X_{0:t}$ given $Y_{0:t}=y_{0:t}$. This posterior is induced by a controlled SDE. Consequently, instead of recursive ensemble updates or iterative sampling, when new observation comes in, our amortized map allows for direct generation of the posterior sample paths and associated uncertainty quantification outputs without re-training.

\subsection{Contribution and main results}

Motivated by analytical formulations of robust nonlinear filtering, we develop a Bayesian inference approach on path space stochastic dynamical systems with noisy and partial measurements. The resulting method combines a stochastic control representation of the posterior path law with an amortized conditional generative model, yielding posterior trajectory and filtering law estimation and uncertainty quantification for new paths. The main contributions are as follows.

\medskip

\begin{itemize}
\item 
We extend the pathwise filtering formulation (\cite{mitter2003variational, van2007filtering}) to the setting of a time-dependent nonlinear observation map $h(t,x)$ with diffusion coefficient $k(t)$. For each observation window $[0,t]$, we identify in Theorem~\ref{thm: pathcontrol} the posterior path law as the law of a controlled diffusion process. Together with a control representation (Theorem~\ref{prop: control}) via a Hopf-Cole transform of the pathwise Zakai equation (Proposition~\ref{prop:path-zakai}), we get a sequence of feedback-controlled SDEs with score-type structure that induce the pathwise filtering law on $[0,t]$ for each $t$.
\item We use the above control representations to construct an observation-conditioned latent SDE for filtering and posterior path inference. The control is parametrized by a neural network and depends on the history of the given observation path, thereby defining an amortized map from noisy trajectories to approximate posterior path laws. We derive a pathwise ELBO in Theorem~\ref{theorem: elbo} which is a variational approximation of the model posterior path law. 
\item 
The transition mechanism or the prior dynamics is not assumed known but learned from data, in contrast to traditional filtering methods. 
Once trained, our model enables fast generation of data-assimilated trajectories for new observations without rerunning a filter online, and it supports uncertainty quantification tasks for nonlinear path functionals (hitting times, occupation times, autocovariances, etc.) rather than only one-time marginals.
We validate the method on stochastic double-well, Lorenz 63, Lorenz 96 and  to MuJoCo Hopper simulator data.
\end{itemize}

\subsection{Related work}
\paragraph{Bayesian filtering methods}

The mathematical foundations are rooted from filtering and smoothing from stochastic analysis and Bayesian inference, see \cite{bain2009fundamentals, reich2015probabilistic} for example.
Computational results \cite{spantini2022coupling, chopin2023computational} relates Bayesian filtering problems to change or transport of measures. In \cite{reich2021fokker}, a Fokker–Planck based interacting particle system preconditioned by empirical covariance matrix is constructed. 
The evolution of the particle system is based on the gradient flow structure that minimizes the Kullback–Leibler divergence between the particle distribution and the desired invariant posterior measures. Recent approaches to data assimilation have introduced Koopman-based methods \cite{chen2025cgkn, chen2025modeling, frion2024neural}, integration of optimal control to particle filters \cite{zhang2023optimal, vargas2023bayesian, yang2013feedback}. 
In particular, feedback particle filter \cite{yang2013feedback} uses the duality between control and filtering and constructs a particle system that interacts through a mean field term. The core difference with our method is that we focus on a more general question of pathwise filtering and path estimation,  through a variational framework instead of direct Monte-Carlo approximation of empirical densities.

\paragraph{Joint parameter and state estimation} 
When the drift and diffusion coefficients in the system, the parameters and states need to be estimated jointly, methods such as extending the state variable to be the pair of state and parameter, smoothing \cite{douc2009forward, briers2010smoothing, boc2013}, Expectation-Maximization (EM) \cite{bocquet2020bayesian, schon2011system, sutter2016variational}, particle Markov Chain Monte-Carlo (MCMC) \cite{andrieu2010pmcmc,fearnhead2008particle} have been used. In particular, \cite{sutter2016variational} develops a variational framework for pathwise nonlinear filtering using mixture of exponential families, and parameter estimation via an EM-type algorithm. The unknown parameters $\theta$ is in the drift and diffusion coefficients in the form of $f(\theta, x)$ where $f(\cdot,\cdot)$ is considered to be given. 

Our method develops a controlled SDE related to the pathwise nonlinear filtering and construct an amortized representation of the posterior path measure. It is flexible in terms of unknown, complicated drift and diffusion coefficients are complicated and the degeneracy of particles is not an issue.

\paragraph{Density propagation in continuous time}

Previous work has also explored the time evolution of posterior density. For example, one approach involves approximating the SPDEs that the (unnormalized) posterior densities via splitting methods \cite{baagmark2023energy, lobbe2023deep}; Alternatively, \cite{pathiraja2021mckean} employs McKean–Vlasov SDEs to characterize the evolution of the filtering measure through interacting particle systems. We share the same focus on a continuous-time measure transport approach but we extend it by learning a generative map of pathwise filtering measure, allowing for amortized inference that bypasses the need for solving (Stochastic) PDEs directly.

\paragraph{Data assimilation via conditional generative models}

Learning-based alternatives to nonlinear filtering have been explored by \cite{le2021deep, qian2022deep}, where Monte Carlo sample paths from a nominal state-space model are used to train a neural network that maps a finite observation window to the corresponding state.
Recently, a growing body of work uses conditional generative models for filtering and data assimilation that aim to represent or sample from conditional filtering or smoothing distributions. 
Generative modeling-based approaches so far utilizes diffusion models \cite{chao2022denoising, bao2023score, rozet2023score, bao2023unified, bao2024ensemble} or Schr\"odinger bridges 
(\cite{shi2022conditional, bao2025ensemble}). In particular, there are typically two time horizons considered, the actual time that the signal and observation dynamics evolve on, and the pseudo-time where the diffusion models or other generative methods evolve on. Through a recursive, prediction-update scheme,
the model serves to learn the prior $p(X_t|Y_t)$, and generative methods to update the posterior $p(X_{t+1}|Y_{t+1})$. This can be integrated well with traditional methods, for example,
integrating score-based diffusion method with particle filter, where the reweighting is realized by training the score function and perform the update of the state at each time by evaluating the reverse-(pseudo) time process.

Our method deviates from this approach, as we consider generative methods on the real time horizon. As noted in \cite{chung2022diffusion}, when designing generative models that evolve in pseudotime, the likelihood term for posterior update is used for every pseudo-time step and thus the posterior sampler is biased except for pseudo-time zero, where the forward state follows the true data distribution. Instead of learning the mapping of an predict-update step, we approximate the posterior SDE given observation paths, avoiding re-training at every update step.

\subsection{Organization of the paper}

The organization of the paper is as follows. 
Section~\ref{sec: setup} gives the standard setup of the stochastic dynamical system that is not directly observable, but with access to a noisy perturbed paths of the underlying true system, we motivate and construct the pathwise filtering and path estimation problem. In Section~\ref{sec: theory}, we derive the pathwise Zakai equation used for our setup and reformulate the problem in a variational way and solve it through stochastic control. 
The optimal control has an explicit dependence on the gradient of the logarithmic of the posterior law. Then, we develop the conditional latent SDE approach that is amortized on observations in Section~\ref{sec: nsde}. One important result is the derivation of the pathwise ELBO loss. Further implementation details of the encoder and the pushforward is introduced in Section~\ref{sec: experiments}.
We perform several experiments on different stochastic datasets (double-well, Lorenz 63, Lorenz 96 systems and MuJoCo Hopper). We illustrate the flexibility and robustness of such methods on estimating path functionals, performing fast inference on time intervals that have not been trained on, missing observations, and comparison with particle-based methods.

\section{Problem set-up}
\label{sec: setup}

In this section, we define the problem setup and motivate our path-dependent approach for data assimilation from both theoretical and practical perspectives.

Under the probability space $(\Omega, \mathcal{F}, \mathbb{P})$, suppose the state of the underlying system follows a general stochastic differential equation (SDE) form
\begin{flalign}
\label{eq: state}
\text{(State, unobservable)} && dX_t = \beta(t, X_t) dt + \sigma(t, X_t)dW_t, \quad t \in [0,T], &&
\end{flalign}
where $ \beta : [0,T] \times \mathbb{R}^{n} \rightarrow \mathbb{R}^n, \sigma : [0,T] \times \mathbb{R}^{n} \rightarrow \mathbb{R}^{n \times p}, \{ W_t \in \mathbb{R}^p \}_{t \in [0,T]} $ denotes a standard Wiener process. We assume the observation $Y \in \mathbb R^m$ follows 
\begin{flalign}
\label{eq: obs}
\text{(Noisy observation)} && d Y_t = h(t,  X_t) dt + k(t) dB_t, \quad t \in [0,T], &&
\end{flalign}
where $h: [0,T] \times \mathbb{R}^{n} \rightarrow \mathbb{R}^m$,  $k:[0,T]\to \mathbb{R}^{m\times m}$ is measurable.  $\{ B_t \in \mathbb{R}^m \}_{t \in [0,T]} $ denotes a standard Wiener process. Define $e(t) := k(t) k(t)^\top$.

\begin{asmp}
\label{asmp:sde}
\begin{enumerate}[label=(\roman*)]
    \item The drift and diffusion $\beta(t,x),\sigma(t,x)$ are measurable in $t$ and $C^{2}$ in $x$, with
all $x$-derivatives up to second order bounded and globally Lipschitz in $x$, uniformly in $t$;
moreover $\beta,\sigma$ satisfy a linear growth condition. The observation function $h(t,x)$ is
of class $C^{1,2}([0,T]\times\mathbb R^{n})$, with $\partial_t h,\ \nabla_x h,\ \nabla_x^{2}h$
bounded uniformly in $(t,x)$.
    \item $a(t,x):=\sigma(t,x)\sigma(t,x)^\top$
  is uniformly elliptic and bounded:
  there exist $0<\lambda\le\Lambda<\infty$ such that
  \[
    \lambda \|v\|^2 \le v^\top a(t,x) v \le \Lambda \|v\|^2
    \quad\text{for all $v\in\mathbb{R}^n$ and $(t,x)\in[0,T]\times\mathbb{R}^n$.}
  \]
\end{enumerate}
\end{asmp}

Under Assumption~\ref{asmp:sde}(i), the system \eqref{eq: state} is well-posed and yields a unique strong solution and that the generalized It\^o formula applies to $h(\cdot,X_\cdot)$. Assumption~\ref{asmp:sde}(ii) will be used later for convexity of Hamiltonian in the derivation of the optimal control.

\subsection{Optimal filter}
Denote $\mathcal{P}(\mathbb R^n)$ as the collection of all Borel probability measures on $\mathbb R^n$. $\hat{X}_t = \E[X_t | \mathcal{F}^Y_t]$ is the best estimate of $X_t$ (almost surely) in the sense that $\E[|X_t - \hat{X}_t|^2] = \inf \{\E[|X_t - x|^2]\}$, for any random variable $x \in L^2(\Omega, \mathcal{F}^Y_t, \mathbb{P})$. More generally, the goal in stochastic filtering is to determine the 'filter', which is the conditional distribution $\pi_t \in \mathcal{P}(\mathbb R^n)$ of the signal $X_t$ at time $t$ given the information accumulated from observing $Y$ throughout the time horizon $[0,t]$. That is,
\begin{equation}
\label{eq:filter}
    \pi_t \varphi = \E[\varphi(X_t) | \mathcal{F}^Y_t],
\end{equation}
where $\varphi: \mathbb R^n \rightarrow \mathbb R$ is a bounded measurable test function, and $\mathcal{F}^Y_t := \sigma(Y_s; 0 \leq s \leq t)$. See \cite{bain2009fundamentals} for a thorough introduction on stochastic filtering.
Essentially, we are dealing with Bayesian inverse problems in every time step of the state process, and we have a sequence (resp. flow) of Bayesian inverse problems in discrete (resp. continuous) time series. 
We aim at estimating the posterior density, and sampling from the posterior density to get a good estimate of the state based on observations. The normalized and un-normalized density has been shown to follow SPDEs called the Kushner equation \cite{doi:10.1137/0302009} and the Zakai equation \cite{zakai1969optimal}, respectively.

\subsection{Motivation for path-based approach}

We aim to learn the SDE dynamics and perform path inference from noisy and transformed observations at the same time. 
Let $\mathcal{C}_T := C([0,T];\mathbb R^{n})$. Recall that the continuous stochastic process $X = \{X(t)\}_{t \in [0,T]}$ on $(\Omega, \mathcal F, \PP)$ is a random variable taking values in $(\mathcal{C}_T, \mathcal{B}(\mathcal{C}_T))$ on $(\Omega, \mathcal F, \PP)$. The pushforward measure $\mathrm{Law}(X) := \PP \circ X^{-1}$ on $(\mathcal{C}_T, \mathcal{B}(\mathcal{C}_T))$ is the path measure of $X$. 
\paragraph{Path dependence}
In many applications, the objects of interest are inherently path dependent. In time series modeling, one needs more than pointwise estimates of the state. For example, probabilities of rare events, hitting or exit times, occupation statistics, autocovariances, or functionals that depend on an entire time window are of keen interest in many dynamical systems especially chaotic and stochastic systems. This can be seen by an Ornstein-Uhlenbeck (OU) process with $X_0 \sim N(0, \frac{\sigma^2}{2\theta})$, and $dX_t = -\theta X_t dt + \sigma dW_t$
leads to the marginal distribution $X_t \sim N(0, \frac{\sigma^2}{2\theta})$ for $t \geq 0$. However, the marginal distributions alone cannot show the mean-reverting property without knowledge of the temporal dynamics.

From a computation cost perspective, conducting Bayesian updates iteratively at each observation time can be expensive and causes degeneracy in particle-based methods.
We can show a relationship between the marginal and the path measure that for $\mu$, $\nu$ taking values in $\mathcal{P}^2(C([0,T]; \mathbb R^n))$, $\sup_{t \in [0,T]} W_{\mathbb R^n, 2}(\mu_t, \nu_t) \leq W_{C([0,T]; \mathbb R^n), 2}(\mu, \nu) $. It is better to optimize the path, which controls the error of each marginal while optimizing marginals where the error can propagate. We can have a similar relationship if the Wasserstein distance is changed to KL divergence, using information processing inequality.

\paragraph{Robustness towards the perturbation in noisy observations}
An interesting and important technical consideration is that by conditioning on an observed path $y$, the mapping from an observation path to the posterior should be well-defined and stable under perturbations of $y$. Real measurement paths typically have finite variation on $[0,T]$, while in the aforementioned setup under the Wiener measure, such paths have measure zero. Thus, it is desirable for $\pi_t(\varphi)$ to be defined
for all continuous paths $y$, not just for those in a set of
probability 1 under the Wiener measure. In the paper, we follow the robust construction of the filter \cite{clark2005robust, crisan2022pathwise} in the sense that the posterior is given by a continuous functional on the path space $C([0,t];\mathbb{R}^{m})$, with respect to the supremum norm.  A slight change in the observed path does not give drastic changes to the filter $\pi_t(\varphi)$.
One 
chooses a version for which $y \mapsto \pi_t \varphi$ extends continuously to all observation paths $y \in C([0,t]; \mathbb{R}^m)$.

\bigskip
Therefore, path-dependent quantities are naturally defined in terms of the conditional law of the latent path given the observations. It is desirable that 1) the filter is not too sensitive towards different discretization of data or perturbations not captured by the observation model. 2) the posterior path measure and the filtering distribution, that is, 
\begin{equation}
\label{eq: pi}
    \Pi_{[0,t]}^y := \Law(X_{[0,t]}|Y_{[0,t]} = y), \qquad \pi_t (\varphi):= \int \varphi(x_t) \Pi_{[0,t]}^y (dx_{[0,t]})
\end{equation}
can be generated accurately without the expensive particle-based updates which causes path degeneracy.
In the following sections, we regard $Y_\cdot$ as a random element in the path space
$C([0,T];\mathbb R^{m})$, with its Borel $\sigma$-algebra.

\section{Amortized variational inference for posterior measure}

\label{sec: theory}

We show in this section that the approximation to the posterior path measure can be solved via a stochastic control problem. This leads to the inference method in Section~\ref{sec: nsde} with principled uncertainty quantification on path space and efficient generation of representative posterior trajectories. The variational Kallianpur-Striebel \cite{kallianpur1969stochastic} can be derived from the Gibbs variational principle (see \cite{van2007filtering, raginsky2024variational}), the posterior measure is achieved by minimizing energy over trajectories. The variational approach to solve nonlinear estimation problems is first presented in \cite{mitter2003variational}. The connection of filtering and control has been discussed in \cite{kim2020optimal, reich2019data}, applied in particle filters \cite{yang2013feedback} and path estimation on exponential families \cite{sutter2016variational}. 
The formulation bridges Gibbs variational inference on path space and controlled diffusions.
However, solving the controlled diffusion directly to generate the posterior path law is computationally undesirable. In this section, we explore the structure of the drift change by a pathwise filtering analysis. We will connect the result with Gibbs variational principle \cite{raginsky2024variational} to develop a neural parametrization for the posterior path laws in the next section.

\subsection{Variational representation posterior path law}
\label{sec: control}

Our main methodological step is to represent, for each fixed horizon $t$, the posterior path law in \eqref{eq: pi} as a controlled diffusion law on $[0,t]$.
Previous work such as \cite{fleming1982optimal} and \cite{mitter2003variational} has proposed this viewpoint and typically connect the optimal control to a gradient of a log-density potential. Under the conditional law of $X$ given $\{Y_t\}_{t \in [0,T]}$, $\{X_t\}_{t \in [0,T]}$ is Markovian for each given path $y$ by \cite[Section 2.2.2]{van2007filtering}. Thus, we discuss controls in feedback form.

Consider a time segment $[s,t] \subset [0,T]$, define the negative log-likelihood ratio
\[
H_{s,t}(X;y)
:=
-\!\int_s^t h(r,X_r)^\top e(r)^{-1}\,dy_r
+ \tfrac12\!\int_s^t \|h(r,X_r)\|^2_{e(r)^{-1}}\,dr ,
\]
where $y\in C^1([0,t];\mathbb R^m)$.
With initial condition $x\in\mathbb R^n$, let $\mu_{s,x}$ be the path law on
$C([s,t];\mathbb R^n)$ induced by the prior dynamics \eqref{eq: state}. Define
$\eta_t^y(s,x)
:=
\E^{\mu_{s,x}}\left[\exp\bigl(-H_{s,t}(X;y)\bigr)\right]$.
Accordingly, define the relative value function
$\widetilde V_t^y(s,x)
:=
-\log \eta_t^y(s,x)$. Thus,
$\widetilde V_t^y(t,x)=0$.
The posterior law is given by Bayes' rule on path space,
\begin{equation}
\label{eq:posterior-segment-forward-fixed}
\frac{d\Pi_{s,x}^{y,t}}{d\mu_{s,x}}(X)
=
\frac{\exp(-H_{s,t}(X;y))}{\eta_t^y(s,x)}.
\end{equation}
Then for every $0\le s\le t$ and $x\in\mathbb R^n$, by the Gibbs variational
principle \cite{raginsky2024variational},
\begin{equation}
\label{eq:gibbs-forward-fixed}
\widetilde V_t^y(s,x)
=
\inf_{\nu\ll \mu_{s,x}}
\left\{
\KL(\nu\|\mu_{s,x})+\E^\nu[H_{s,t}(X;y)]
\right\},
\end{equation}
and the unique minimizer is $\nu=\Pi_{s,x}^{y,t}$.

\begin{thm}[Fixed-horizon forward variational representation of the posterior path law]
\label{thm: pathcontrol}
Fix $t\in(0,T]$ and an observation path $y\in C^1([0,t];\mathbb R^m)$.
Under Assumption~\ref{asmp:sde}-\ref{asmp: novikov}, assume
$\widetilde V_t^y\in C^{1,2}([0,t)\times\mathbb R^n)\cap
C([0,t]\times\mathbb R^n)$. Let $\QQ_{s,x}^{\alpha}$ be the path law induced by
\begin{equation}
\label{eq:controlled-sde-forward-fixed}
dX_r
=
\bigl(\beta(r,X_r)+\sigma(r,X_r)\alpha_r\bigr)\,dr
+
\sigma(r,X_r)\,dW_r^{\alpha},
\qquad r\in[s,t],
\qquad X_s=x,
\end{equation}
where $\alpha$ is progressively measurable and
$\E^{\QQ_{s,x}^{\alpha}}\!\left[\int_s^t \|\alpha_r\|^2\,dr\right]<\infty$.
The unique optimal feedback is
\begin{equation}
\label{eq:alpha-star-forward-fixed}
\alpha_t^{y,\star}(s,x)
=
-\sigma(s,x)^\top \nabla_x \widetilde V_t^y(s,x).
\end{equation}
The law of the controlled diffusion \eqref{eq:controlled-sde-forward-fixed}
with optimal control \eqref{eq:alpha-star-forward-fixed} is $\Pi_{s,x}^{y,t}$.
\end{thm}

\begin{proof}
By Girsanov's theorem, $\KL(\QQ_{s,x}^{\alpha}\|\mu_{s,x})
=
\frac12\E^{\QQ_{s,x}^{\alpha}}\!\left[\int_s^t \|\alpha_r\|^2\,dr\right].$ Substituting this into the Gibbs formula \eqref{eq:gibbs-forward-fixed} yields
\begin{equation}
\label{eq:control-forward-fixed}
\widetilde V_t^y(s,x)
=
\inf_{\alpha}
\E_{s,x}^{\QQ^\alpha}\left[
\frac12\int_s^t \|\alpha_r\|^2\,dr
+
H_{s,t}(X;y)
\right].
\end{equation}
Let
$L(r,x)
=
\frac12\|h(r,x)\|_{e(r)^{-1}}^2
-\dot y_r^\top e(r)^{-1}h(r,x)$.
The controlled generator acts on smooth $\varphi$ as
\[
\mathcal L_r^\alpha \varphi(x)
=
(\beta(r,x)+\sigma(r,x)\alpha)\cdot \nabla_x \varphi(x)
+\frac12\mathrm{tr}\!\left(a(r,x)\nabla_x^2\varphi(x)\right).
\]
Applying the dynamic programming principle on $[s,s+\delta]$, and letting $\delta\downarrow0$, yields
\[
-\partial_s \widetilde V_t^y
=
\frac12\mathrm{tr}(a\nabla_x^2\widetilde V_t^y)
+
\inf_{\alpha}
\left\{
(\beta+\sigma\alpha)\cdot \nabla_x \widetilde V_t^y
+\frac12\|\alpha\|^2
+L
\right\}.
\]
The unique minimizer is \eqref{eq:alpha-star-forward-fixed}. 

Now, apply It\^o's formula to $r\mapsto \widetilde V_t^y(r,X_r)$, we obtain
\[
d\widetilde V_t^y(r,X_r)
=
\left(
\frac12\|\alpha_r^\star\|^2
-
L(r,X_r)
\right)\,dr
-
(\alpha_r^\star)^\top dW_r.
\]
Since $\widetilde V_t^y(t,\cdot)=0$, by integrating from $s$ to $t$ yields
\[
-H_{s,t}(X;y) +\widetilde V_t^y(s,x)
=
\int_s^t (\alpha_r^\star)^\top\,dW_r
-\frac12\int_s^t \|\alpha_r^\star\|^2\,dr.
\]
Thus,
\[
\frac{\exp(-H_{s,t}(X;y))}{\eta_t^y(s,x)}
=
\exp\left(
\int_s^t (\alpha_r^\star)^\top\,dW_r
-\frac12\int_s^t \|\alpha_r^\star\|^2\,dr
\right)
=
\frac{d\QQ_{s,x}^{\alpha^\star}}{d\mu_{s,x}}.
\]
Comparing with \eqref{eq:posterior-segment-forward-fixed} shows that $\QQ_{s,x}^{\alpha^\star}=\Pi_{s,x}^{y,t}$.
\end{proof}

So far, we have identified the posterior path law on $[0,t]$ as the law of an optimally controlled process, with feedback control that is causal in $y$. The variational formulation and the explicit controlled SDE leads to a practical algorithm in the next section. We next look further into the filtering density.

\subsection{Pathwise Zakai equation}

Building on the robust filter viewpoint, for $\mathbb{P}_Y-a.e.$ observation path $y$,
we first derive a parabolic PDE that governs the evolution of an unormalized conditional density as a functional of a deterministic observation path in the following proposition.
A standard route is through the unnormalized (Zakai) functional $\rho_t^y$ and the Kallianpur-Striebel normalization,
\begin{equation}
\label{eq: q-orig}
    \pi_t^y(\varphi)=\frac{\rho_t^y(\varphi)}{\rho_t^y(1)},\qquad 
\rho_t^y(\varphi)=\int_{\mathbb{R}^n} \varphi(x)\,r_t^y(x)\,dx,
\end{equation}
where $r_t^y$ is the unnormalized pathwise conditional density.

We follow the change of measure approach in \cite[Theorem 5.12]{bain2009fundamentals} and pathwise diffusion filtering result in \cite{clark2005robust, van2007filtering} to derive the pathwise Zakai equation for nonlinear observation model with time dependent diffusion coefficient $k(t)$ as defined in \eqref{eq: obs}. We show in the following proposition that the unnormalized density follows a parabolic PDE. 

\begin{asmp}[Novikov's condition]
\label{asmp: novikov}
$k(t)$ is an invertible, measurable function for all $t \in [0,T]$. In addition, it holds that
  \[
    \E\exp\left(
      \frac{1}{2}\int_0^T
        \|k(t)^{-1}h(t,X_t)\|^2\,dt
    \right) < \infty;\]   
and for every deterministic path
  $y\in C([0,T];\mathbb{R}^{m})$,
\[
\E\exp\left(
      \frac{1}{2}\int_0^T
        \left\|
         \left((\nabla_x h)(t,X_t)\,\sigma(t,X_t)\right)^\top e(t)^{-1} y_t
        \right\|^2\,dt
    \right) < \infty.
\]
\end{asmp}

\begin{prop}
\label{prop:path-zakai}

Under Assumption~\ref{asmp:sde}-\ref{asmp: novikov}, for observation $y \in C([0,T];\mathbb{R}^{m})$ a.e., let $q_t^y = \exp \{ - h(t,x)^\top e(t)^{-1}y_t\} r_t^y$, where $r_t^y$ is the unnormalized posterior density in \eqref{eq: q-orig}. It follows
\begin{equation}
\label{eq: pathwisezakai}
\partial_t q^y_t = - \nabla \cdot \left( \beta^y(t, x) q^y_t \right) + \frac{1}{2} \sum_{i,j=1}^n \partial_{x_i x_j}(a_{ij}(t, x) q^y_t) + q^y_t G(t, x; y),
\end{equation}
where $a(t,x) = \sigma(t,x) \sigma(t,x)^\top$,
$\beta^y (t,x):=\beta(t,x)-a(t,x)(\nabla_x h(t,x))^\top e^{-1}(t) y_t$,
\begin{equation}
\label{eq: G}
\begin{aligned}
   & G(t,x;y) =
-\frac12\|h(t,x)\|_{e(t)^{-1}}^2
-y_t^\top e(t)^{-1}\left(\partial_t h(t,x)+\mathcal L_t h(t,x)\right)
\\
& \qquad \qquad \qquad  -y_t^\top \dot e(t)^{-1} h(t,x) +\frac12\Big\|((\nabla_x h)(t,x)\sigma(t,x))^\top e(t)^{-1}y_t\Big\|^2,
\end{aligned}
\end{equation}
where $\dot e(t)^{-1} := \frac{d}{dt}(e(t)^{-1})$ and
\begin{equation}
\label{eq: generator}
    \mathcal{L}_t f(x) = \sum_{i=1}^{n} \beta_i(t,x) \partial_{x_i} f(x) + \frac{1}{2} \sum_{i,j=1}^{n} a_{ij}(t,x) \partial_{x_i x_j} f(x)
.
\end{equation} 
\end{prop}

\begin{proof}
First, we perform a change of measure such that \eqref{eq: obs} depends only on time and the new Brownian motion under the new measure $\mathbb{Q}$. In this way, \eqref{eq: obs} is independent of $X$.
Write $\|u\|_{e(t)^{-1}}^2:=u^\top e(t)^{-1}u$.
Define the Radon-Nikodym derivative by
\[
\frac{d\QQ}{d\PP}\Big|_{\mathcal{F}_t}
= \exp \left(-\int_0^t h(s,X_s)^\top k^{-\top}(s) dB_s
-\frac{1}{2}\int_0^t \|h(s,X_s)\|^2_{e(s)^{-1}} ds\right).
\]
Then, by Assumption~\ref{asmp: novikov} and Girsanov theorem, the process
\[
\widetilde B_t := B_t + \int_0^t \theta_s\,ds, \quad \theta_t := k^{-1} h(t,X_t)\in\mathbb{R}^m
\]
is an $\mathbb{R}^m$–Brownian motion under $\QQ$, and the observation becomes driftless $dY_t = k(t)\,d\widetilde B_t$.
Because the density depends only on $B$ (and on $Y$), $W$ remains a Brownian motion under $\QQ$ and is independent of $\widetilde B$. Hence the process $X$ driven by $W$, and the process $Y$ as a deterministic transform of $\widetilde B$ are independent under $\QQ$. Moreover, the law of $(X,\widetilde B)$ under $\PP$ equals the law of $(X,Y)$ under $\QQ$.
As $h(t, X_t)$ is a semimartingale, we have
\begin{equation}
    \label{eq:pathh}
    \begin{aligned}
 &     \int_0^t h(s,X_s)^\top e(s)^{-1} dY_s \\
& =
h(t,X_t)^\top e(t)^{-1}Y_t
-
\int_0^t Y_s^\top e(s)^{-1}\,dh(s,X_s)
-
\int_0^t Y_s^\top \dot e(s)^{-1}h(s,X_s) ds,
    \end{aligned}
    \end{equation}
    $\mathbb{Q}$-almost surely. {The cross-variation term $\langle h(X), Y \rangle \equiv 0$ in \eqref{eq:pathh} under $\mathbb{Q}$, as $W$ and $\widetilde{B}$ are independent.} This is a consequence from \cite[theorem 5.12]{bain2009fundamentals}. In particular, $ \int_0^T Y_s dh(X_s)$ is well-defined as $Y_s$ is $\{\mathcal{F}_t\}$-adapted. By Kallianpur-Striebel formula, 
    \[
    \begin{aligned}
    \rho_t(f) 
        =& \E^{{\mathbb{Q}}}\left[ f({X}_t) \exp\left( \int_0^t h(s, {X}_s){e(s)^{-1}} dY_s  - \frac{1}{2} \int_0^t \|h(s, X_s)\|^2_{e(s)^{-1}} ds\right)\right].
    \end{aligned}
\]
We now perform the second change of measure through the Girsanov theorem. Apply generalized It\^o's formula, we get
\[
dh(s,X_s) = \partial_s h(s, X_s) \md s + \mathcal{L}_s h(s, X_s) \md s +  (\nabla_x h)(s,X_s) \sigma(s, X_s) dW_s.
\]
{Thus, define
$\mu_s(x;y):=\left((\nabla_x h)(s,x)\,\sigma(s,x)\right)^\top e(s)^{-1} y_s\in\mathbb R^p$,
and set the Radon-Nikodym derivative
\begin{equation}
\label{eq: rn2}
    \frac{d\widetilde{\mathbb P}}{d\mathbb Q}\Bigg|_{\mathcal F_t}
=\exp\left(-\int_0^t \mu_s(X_s;y)^\top dW_s-\frac12\int_0^t \|\mu_s(X_s;y)\|^2 ds\right),
\end{equation}
so that $\widetilde W_t:=W_t+\int_0^t \mu_s(X_s;y)\,ds$ is a Brownian motion under $\widetilde{\mathbb P}$. Under the new probability measure $\widetilde{\mathbb{P}}$, the state process $\widetilde{X}_t$ has a modified drift, i.e.,
\[
d\widetilde{X}_t =\left[\beta(t,\widetilde{X}_t )-a(t,\widetilde{X}_t )(\nabla_x h)(t,\widetilde{X}_t )^\top e(t)^{-1} y_t\right]dt+\sigma(t,\widetilde{X}_t )\, d\widetilde W_t, \qquad \widetilde{X}_0 = X_0.
\]} 
Thus, for $f \in C^\infty_c(\mathbb{R}^n)$,
\begin{equation}
\label{eq: tau}
    \rho_t^y(f) = \mathbb{E}^{\widetilde{P}}\left[f(\widetilde{X}_t) \exp \left \{h(t,\widetilde X_t)^\top e(t)^{-1} y_t\right\}\exp\left(\int_0^t G(s, \widetilde{X}_s; y) ds\right)\right].
\end{equation}
Define for every $t\in[0,T]$,
\begin{equation}
\label{eq:q-def}
\begin{aligned}
  \int_{\mathbb R^n}f(x) q_t^y(x)\,dx
    :=
    \mathbb E^{\widetilde{\mathbb P}} \left[
       f(\widetilde X_t)\exp \left(\int_0^t G(s,\widetilde X_s;y)\,ds\right)
    \right],
    \qquad f\in C_c^\infty(\mathbb R^n).
  \end{aligned}
\end{equation}
Assume that, for
the given observation path $y$, 
$q^y\in C^{1,2}([0,T]\times\mathbb{R}^n)$. By \eqref{eq: q-orig}, we get $q_t^y = \exp \{ -h(t,x)^\top e(t)^{-1}y_t\} r_t^y$.

Apply It\^o’s formula under $\widetilde{\PP}$ to $Z_t := f(\widetilde X_t)\exp \left(\int_0^t G(s,\widetilde X_s;y)\,ds\right)$, we get
\[
\begin{aligned}
dZ_t
&= \exp \left(\int_0^t G(s,\widetilde X_s;y)\,ds\right)
   \left[
      \left(\mathcal{L}_t^y f + G(t,\cdot;y)\,f\right)(\widetilde X_t)\,dt
      + (\nabla_x f(\widetilde X_t))^\top\sigma(t,\widetilde X_t)\,d\widetilde W_t
   \right],
\end{aligned}
\]
where \[
(\mathcal L_t^y) f(x):=\sum_{i=1}^n \beta_i^y(t,x)\,\partial_{x_i} f+\frac12\sum_{i,j=1}^n a_{ij}(t,x)\,\partial_{x_i x_j}^2 f.
\]
Taking expectations under $\widetilde{\PP}$ and using that the stochastic
integral has mean zero yields
\[
  \frac{d}{dt}\,
  \mathbb{E}^{\widetilde{\PP}} \left[
    f(\widetilde X_t)\,
    e^{\int_0^t G(s,\widetilde X_s;y)\,ds}
  \right]
  = \mathbb{E}^{\widetilde{\PP}} \left[
      \left(\mathcal{L}_t^y f + G(t,\cdot;y)\,f\right)(\widetilde X_t)\,
      e^{\int_0^t G(s,\widetilde X_s;y)\,ds}
    \right].
\]
By \eqref{eq: tau}, the left-hand side is
$\frac{d}{dt}\int f(x)\,q_t^y(x)\,dx$, hence, with $q^y\in C^{1,2}$,
\begin{equation}
\label{eq:weak-evolution}
  \frac{d}{dt}\int_{\mathbb{R}^d} f(x)\,q_t^y(x)\, \md x
  = \int_{\mathbb{R}^d} f(x)\,\partial_t q_t^y(x)\, \md x =  \int_{\mathbb{R}^d}
      \left(\mathcal{L}_t^y f(x) + G(t,x;y)\,f(x)\right)\,
      q_t^y(x)\, \md x
\end{equation}
for all $f\in C_c^\infty(\mathbb{R}^d)$.
Integrating by parts in $x$, we conclude the unnormalized filtering density \eqref{eq: pathwisezakai}, for a.e.\ path $y$.
\end{proof}

\begin{definition}[Admissible controls]
The feedback control $u:[0,T]\times\mathbb R^n\to\mathbb R^n$ is called admissible if it is continuous, satisfies a uniform Lipschitz condition on every compact set for $t \in [0,T]$, and the linear growth condition. That is, for every compact set $K\subset\mathbb R^n$, there exists $L_K<\infty$ with
    \[
    |u(t,x)-u(t,x')|\le L_K |x-x'|,
    \qquad t\in[0,T],\ x,x'\in K;
    \]and there exists $C_u < \infty$ such that $|u(t,x)|\le C_u(1+|x|)$ for all $(t,x)\in[0,T]\times\mathbb{R}^n$.
We denote admissible controls as $u \in U$.
\end{definition}
The above defines the feedback control that ensures the controlled SDE is well-posed and the running cost has at most quadratic growth.
{Assume $q^y_0 > 0$ on the support of interest so that $V$ is well defined.} We now show that we can construct a HJB equation and an associated control problem from the pathwise Zakai equation. It gives the optimal feedback explicitly in the form of the score function of pathwise posterior density. 

\begin{thm}
\label{prop: control}
Given a continuous path $y := \{Y_t\}_{t \in [0,T]}$, take $V(t,x) = -\log q^y_t(x) \in C^{1,2}([0,T] \times \mathbb{R}^n; \mathbb{R})$, where the unnormalized density $q_t$ solves \eqref{eq: pathwisezakai} in the weak sense.
Then $V$ satisfies a nonlinear parabolic equation
\begin{equation}
\label{eq: parab}
    \partial_t V_t= \frac{1}{2} \mathrm{tr}\left(a(t,x)\nabla_x^2 V_t\right) + H(t, x, \nabla_x V_t),\quad V(0,x) = -\log q^y_0(x).
\end{equation}
where $H(t, x, p)= g^y(t,x) \cdot p - \frac{1}{2} \|\sigma(t,x)^\top p\|^2 - \widehat{G}(t, x; y)$, with $g^y(t,x) := (\nabla \cdot a)(t, x)-\beta^y(t, x)$, $\widehat{G}(t, x; y) = G(t, x; y) - \Theta(t, x; y)$,
$\Theta(t, x; y) := \nabla \cdot \beta^y(t, x) - \frac{1}{2} \sum_{i,j =1}^n \partial_{x_i x_j} a_{ij}(t, x)$, $G(t, x; y)$ is defined in \eqref{eq: G}, and $\mathcal{L}(\cdot)$ defined in \eqref{eq: generator}. 

Define the running cost 
\begin{equation}
\label{eq:running-cost-clean}
\ell(t,x,u)
:=
\frac12\,\langle a(t,x)^{-1}(u-g^y(t,x)),u-g^y(t,x)\rangle
-\widehat G(t,x;y),
\end{equation}
the minimizer is
\begin{equation}
\label{eq:optcontrol}
\tilde{u}^\star(t,x)=g^y(t,x)-a(t,x)\nabla_xV(t,x).
\end{equation}

For $t\in(0,T]$ and initial state $x\in\mathbb R^n$, we consider a time-reversed process with feedback form admissible control $\tilde{u}(\cdot) \in U$,
    \begin{equation}
    \label{eq:controlled-sde-tau}
    d\widetilde{X}_s^u=\tilde{u}(t-s,\widetilde{X}_s^u)\,ds+\sigma(t-s,\widetilde{X}_s^u)\,dW_s,
    \qquad
    \widetilde{X}_0^u=x, \quad s \in [0,t],
    \end{equation}
and with $\phi(x):=-\log q^y_0(x)$, the finite-horizon value function is
    \begin{equation}
    \label{eq:value-tau}
    \mathcal V(t,x)
    =
    \inf_{u\in U}
    \mathbb E_x \left[
    \int_0^t \ell\left(t-s,\widetilde{X}_s^u,u(t-s,\widetilde{X}_s^u)\right)\,ds+\phi(\widetilde{X}_t^u)
    \right].
    \end{equation}
Then, with the unique optimal control \eqref{eq:optcontrol}, $V=\mathcal V$, and \eqref{eq: parab} is the HJB equation.
\end{thm}

\begin{proof}
Since $q^y \in C^{1,2}$ and is strictly positive, we perform Hopf-Cole transform and let
$V(t,x):=-\log q^y_t(x)$, where $V\in C^{1,2}([0,T]\times\mathbb{R}^n)$. We compute the following derivatives on $q^y=e^{-V}$,
\begin{equation}
    \label{eq:qs}
    \partial_t q^y=-q^y\,\partial_tV,
\qquad
\nabla_x q^y=-q^y\,\nabla_xV,
\qquad
\partial_{x_ix_j}q^y
=
q^y\left(\partial_{x_i}V\,\partial_{x_j}V-\partial_{x_ix_j}V\right).
\end{equation}
Expand $-\nabla\cdot(\beta^y q^y)$ and $\frac12\sum_{i,j=1}^n \partial_{x_ix_j}(a_{ij}q^y)$ in \eqref{eq: pathwisezakai}, we get
\[
-\nabla \cdot(\beta^y q)
=
-\sum_{i=1}^n \partial_{x_i}(\beta_i^y q)
=
-\sum_{i=1}^n (\partial_{x_i}\beta_i^y)\,q
-\sum_{i=1}^n \beta_i^y\,\partial_{x_i}q,
\]
\[
\frac12\sum_{i,j=1}^n \partial_{x_i x_j}(a_{ij}q)
=
\frac12\sum_{i,j=1}^n
\left[
(\partial_{x_i x_j}a_{ij})\,q
+
(\partial_{x_j}a_{ij})\,\partial_{x_i}q
+
(\partial_{x_i}a_{ij})\,\partial_{x_j}q
+
a_{ij}\,\partial_{x_i x_j}q
\right].
\]
Further, by substituting the corresponding terms with \eqref{eq:qs}, we derive that $V$ solves the viscous Hamilton-Jacobi equation with \eqref{eq: parab}.
For any fixed $(t,x,p)$, the map 
\[
F(u) := u \cdot p + \ell(t,x,u) = u \cdot p + \frac12\|u-g^y(t,x)\|_{a(t,x)^{-1}}^2 - \widehat G(t,x;y)
\] 
is strictly convex in $u$ by the uniform ellipticity of $a$. Thus, by minimizing the terms on $u$, we get that the unique minimizer is $\tilde{u}^\star = g^y(t,x) - a(t,x)p$ for the Hamiltonian $H(t,x, p) = \inf_u F(u)$.

In \eqref{eq:running-cost-clean}, by the growth assumptions Assumption~\ref{asmp:sde} on $\beta,\sigma$ and $u\in U$, there exists $C$ such that for all $(t,x,u) \in [0,T) \times \mathbb{R}^n \times U$,
$
|\ell(t,x,u)|
\le C\left(1+|x|^2\right)$,
$|V(t,x)|\le C\left(1+|x|^2\right)$.
Thus, $V$ and $\ell$ have at most quadratic growth and $\E_{t,x}\left[\sup_{s\in[t,T]}|\widetilde{X}_s^u|^2\right]<\infty$.
Hence, for every admissible $u$, the controlled state process is well-posed.
Next, fix an admissible control $u$. By the regularity in Assumption~\ref{asmp:sde}, and $V(t,x):=-\log q^y_t(x)$, we apply It\^o's formula to $V(t-s,X_s^u)$ for $0\le s\le t$, and it follows
\[
\begin{aligned}
dV(t-s,\widetilde{X}_s^u)
&=
\Big(
-\partial_t V(t-s,\widetilde{X}_s^u)
+
\tilde{u}(t-s,\widetilde{X}_s^u)\cdot \nabla_xV(t-s,\widetilde{X}_s^u) \\
&\qquad +
\frac12\operatorname{Tr} \big(a(t-s,\widetilde{X}_s^u)\nabla_x^2V(t-s,\widetilde{X}_s^u)\big)
\Big)\,ds
+
(\nabla_xV)^\top \sigma(t-s,\widetilde{X}_s^u)\,dW_s.
\end{aligned}
\]
By \eqref{eq: parab},
\[
-\partial_t V(t-s, x)
+
\tilde{u}(t-s,x)\cdot \nabla_xV(t-s, x)
+
\frac12\operatorname{Tr}(a(t-s, x) \nabla_x^2V(t-s, x))
+
\ell(t-s,x,u)
\ge 0.
\]
Integrating from $0$ to $t$, taking expectations, and using the martingale property of the stochastic integral, we obtain
\[
V(t,x)\le
\mathbb E_x \left[
\int_0^t \ell\left(t-s,\widetilde{X}_s^u,\tilde{u}(t-s,\widetilde{X}_s^u)\right)\,ds+\phi(\widetilde{X}_t^u)
\right].
\]
If $\tilde{u}=\tilde{u}^\star$, then the pointwise Hamiltonian inequality becomes equality, hence equality holds above, $V=\mathcal V$, and $u^\star$ is optimal. Uniqueness follows from strict convexity of the Hamiltonian in $u$.
\end{proof}

In summary, Theorem~\ref{thm: pathcontrol}
shows that the posterior path law can be represented as a controlled
diffusion path law, and Proposition~\ref{prop:path-zakai} and Theorem~\ref{prop: control} identifies the
corresponding filtering correction as a score-type feedback given by the solution of the pathwise Zakai PDE. The two theorems above are linked through a forward-backward factorization of the unnormalized posterior
density. Indeed, by the Markov property under the reference law and the multiplicative structure of the
likelihood functional, $\exp(-H_{0,t}(X;y))=\exp(-H_{0,s}(X;y))\exp(-H_{s,t}(X;y))$;
conditioning on $X_s=x$ and taking the reference expectation of the second factor
gives $\eta_t^y(s,x)$. Hence, up to an $x$-independent constant, $q_s^{y_{[0,t]}}(x) \propto q_s^{y_{[0,s]}}(x) \eta_t^y(s,x)$. Together, these results motivate the inference as a sequence of subproblems:ranging $t$ over the time grid yields the family
$\{\Pi^{y,t_k}_0\}_k$, one posterior path law per horizon.

\section{Learning path estimation through conditional neural SDEs}
\label{sec: nsde}

We present a practical algorithm for learning an amortized approximation of the pathwise filter and the posterior path law solved in the last section. 
{As shown in Theorem~\ref{thm: pathcontrol}, the pathwise filtering problem can be interpreted as constructing a feedback drift correction that transports a reference diffusion toward the conditional law determined by the observation path $y$. We thus consider a parametrized family of conditional laws given observation paths $y$ to approximate $\Pi^y$.}
Note that in traditional filtering methods, albeit the non-observability of the signal process, the signal dynamics or the transition function is given. In contrast, we simultaneously learn a dynamical generative model for the signal evolution and variational approximation to its posterior path law given noisy observations. Throughout this section, we discretize at times $0=t_0<\dots<t_N=T$ and write
the signal path on the grid as $X_{0:N}:=(X_{t_k})_{k=0}^N$.

\subsection{Conditional latent SDEs}

We construct a dynamic generative approach based on neural SDEs \cite{kidger2022neural} and the variational structure in \cite{li2020scalable}. The idea is to approximate the evolution of the posterior state $X_t^y$ with a conditional stochastic map dependent on the observation path $y$, and we represent such state through a latent process $Z$, which evolves in the latent space $\mathcal Z=\mathbb R^{d_\ell}$. 
We learn a conditional latent path measure $\PP_\theta^{ y}$ on $\Omega := C([0,T];\mathcal{Z})$, and a decoder $\mathcal{D}_\theta : \mathcal{Z} \to \RR^n$ such that the decoded law $\widehat{\Pi}^{y}_\theta := (\mathcal{D}_\theta)_\sharp \mathbb{P}_\theta^{y}$ is an amortized approximation to the filtering posterior law $\Pi^y$ defined in Section~\ref{sec: theory}.
During training, paired sequences $(x_{0:N},y_{0:N})$ are available. 
We introduce an auxiliary
SDE whose solution induces a {conditional path measure}
$\mathbb{Q}_\phi^{x, y}$ on $\Omega$. 
At test time, given a new observation path $y$, we sample from $\widehat{\Pi}^{y}_\theta$.

\paragraph{Encoder and initial latent states}
We model the training-time initial latent law by a conditional density $Z_0\sim q_{\phi,0}(\cdot\mid x_{0:N},y_{0:N})$,
parameterized by an encoder network ${\varepsilon}_\phi$. In particular, we encode the initial condition $z_0$ of the approximate posterior process \eqref{eq:postz} as
\[
q_\phi(z_0 \mid x_{0:N},y_{0:N})
=
\mathcal{N} \left(z_0;\,\mu_\phi(x_{0:N},y_{0:N}),\,\Sigma_\phi(x_{0:N},y_{0:N})\right),
\]
using the reparameterization trick $z_0 = \mu_\phi(x, y) + \Sigma_\phi(x, y)^{1/2} \cdot \epsilon_z$, $\epsilon_z \sim \mathcal{N}(0, I)$. 
So ${\varepsilon}_\phi$ maps the entire observation path to a distribution over the initial latent state $z_0$.

\paragraph{Latent dynamics}

We define the auxiliary latent process as follows.
Let $Z_t$ solve
\begin{equation}\label{eq:postz}
dZ_t
=
f_{\phi}(t,Z_t;x_{0:N}, y_{[0,t]})\,dt
+ g_{\theta}(t,Z_t)\,dW_t,
\qquad Z_0\sim q_\phi(\cdot\mid x,y),
\end{equation}
Define $\mathbb{Q}_\phi^{x,y}$ as the path law on $C([0,T];\mathcal{Z})$ induced by \eqref{eq:postz}.
For each $t$, let $q_{\phi,t}^y$ be the time-$t$ marginal of $\mathbb{Q}_\phi^{x,y}$.
Conditional on a Brownian realization $\omega$, strong well-posedness of \eqref{eq:postz} yields a
measurable solution map $\Phi_{0,t}^{\theta,\phi}(\omega):\mathcal Z\to\mathcal Z$ with
$Z_t=\Phi_{0,t}^{\theta,\phi}(\omega)(Z_0)$; With $Z_0$ independent of $W$, the conditional law of $Z_t$ is the pushforward averaging over the Brownian motion, that is,
$q_{\phi,t}^{x,y}
=\mathbb E_W\!\Big[\left(\Phi_{0,t}^{\theta,\phi}(W)\right)_\sharp\,q_{\phi,0}^{x,y}\Big]$.

\paragraph{Generator}

In addition to the SDE system \eqref{eq:postz} on the latent space, we construct another latent process serves as the generator for the (filtering) posterior path measure.
\begin{align}
    \label{eq:prior-sde}
    dZ_t &= w_\theta(t, Z_t; y_{[0,t]}) dt + g_{\theta}(t,Z_t)\,u_\phi(t,Z_t;y_{[0,t]}) dt + g_\theta(t, Z_t) dW_t, & Z_0 &\sim p_\theta(z_0)
\end{align}
where the parametric family $u_\phi(t,Z_t;y_{[0,t]})$ plays the role of an amortized feedback control in Theorem~\ref{prop: control}, and it depends on the historical observation path $y$ up to $t$. Thus, the parametric latent generative model \eqref{eq:prior-sde} induces the latent path measure $\mathbb{P}^y_\theta$ on $C([0,T];\mathcal{Z})$. For brevity, we denote $\mu_{\theta}(t,z;y_{[0,t]})
:=
w_\theta(t,z;y_{[0,t]})
+ g_{\theta}(t,Z_t)\,u_\phi(t,Z_t;y_{[0,t]})$. {In a setup where the prior drift in \eqref{eq: state} is known, we can incorporate the given structure in $w_\theta$.}

A decoder $\mathcal{D}_\theta:\mathcal{Z}\to\mathbb{R}^n$ maps latent states to signal states, $x_t=\mathcal{D}_\theta(z_t)$. The decoded pushforward $(\mathcal D_\theta)_\sharp \mathbb P_\theta^y$ serves as our approximation to $\Pi^y$. We emphasize that the prior dynamics \eqref{eq: state} and the prior law $\mu$ are not learned separately. Under standard regularity conditions for the well-posedness of SDEs, nondegenerate diffusions, and Novikov's condition, the measures satisfy the absolute continuity
$\mathbb{Q}_\phi^{x, y}\ll \mathbb{P}_\theta^{y}$, and $\mathrm{KL}(\mathbb{Q}_\phi^{x, y}\|\mathbb{P}_\theta^{y})$ is computable. Indeed, let ${\alpha}_{\phi,\theta} = g_\theta^{-1}(t, Z_t) \left( \mu_\theta(t, Z_t; y) - f_{\phi}(t, Z_t; x,y) \right)$. With $\E \left[\int_0^T \frac{1}{2}\|{\alpha}_{\phi,\theta}(t,Z_t;y)\|^2 dt \right] < \infty$, it holds
\begin{equation}
\label{eq:variation-rn}
\frac{d \mathbb{Q}_\phi^{x,y}}{ d \mathbb{P}^y_\theta}(Z)
=
\exp \left(
\int_0^T {\alpha}_{\phi,\theta}(t,Z_t; x, y)^\top \, dW_t
-\frac{1}{2} \int_0^T \|{\alpha}_{\phi,\theta}(t,Z_t;x, y)\|^2 dt
\right)\,
\frac{dq_\phi(Z_0\mid x, y)}{dp_\theta(Z_0)}.
\end{equation}
This will appear in the training objective as we specify later in this section.

\subsection{The pathwise evidence lower bound (ELBO)}

We next derive the training objective over the parameters $(\theta,\phi)$. This is an extension based on the \cite{tzen2019neural, li2020scalable} to suit our filtering purpose.
 At test time, given a new observation path $y$, we generate posterior trajectory samples from the learned conditional generative dynamics \eqref{eq:prior-sde} that yields
${\mathcal{D}_\theta}_\sharp\mathbb{P}_\theta^{\,y}$.
Consequently, maximizing $\mathcal L(\theta,\phi)$ over $\phi$ (with $\theta$ fixed) is equivalent to minimizing the reverse KL to the model posterior. We first impose the usual variational
factorization as follows.
\begin{asmp}[Variational factorization]\label{asmp:var_prior}
Conditionally on the latent path $z_{0:N}$, the states factorize as
\begin{equation}\label{eq:asmplatent}
  p_\theta(x_{0:N}\mid z_{0:N})
  = \prod_{i=0}^N p_\theta(x_i\mid z_{t_i}).
\end{equation}
\end{asmp}

\begin{thm}[Pathwise ELBO]
\label{theorem: elbo}
Under Assumption~\ref{asmp:var_prior}, fix an observation path $y$ and assume $\mathbb{Q}_\phi^{x, y}\ll \mathbb{P}_\theta^{y}$. The joint log-density of the signal and observation is lower-bounded by the functional $\mathcal{L}(\theta, \phi; x, y)$, that is,
\begin{equation}
    \label{eq: condelbo}
\log p_\theta(x_{0:N}|y_{0:N})
\ge
\mathbb E^{\mathbb Q_\phi^{x, y}} \left[
\sum_{i=1}^N \log p_\theta(x_i\mid Z_{t_i})
\right]
-\mathrm{KL} \left(\mathbb Q_\phi^{x, y}\,\|\,\mathbb P_\theta^y\right) =: \mathcal{L}(\theta, \phi; x_{0:N},y_{0:N}).
\end{equation}
Furthermore, for fixed $\theta$ and data, ELBO training over $\phi$ is exactly variational approximation of the model posterior
path law,
    \begin{equation}\label{eq:phi_opt_is_KLproj}
\arg\max_{\phi}\mathcal{L}(\theta,\phi; x,y)
\;=\;
\arg\min_{\phi}\mathrm{KL} \left(Q_\phi^{x,y}\,\|\,\Pi_\theta^{x, y}\right),
\end{equation}
where the model posterior path law
$\Pi_\theta^{x,y}(dz) := \PP_\theta^{y}(dz\mid x_{0:N})
\;\propto\;
p_\theta(x_{0:N} \mid z)\,\mathbb P^y_\theta(dz)$.
\end{thm}

\begin{proof}
Under the shared diffusion and suitable regularity, the path law $\mathbb{Q}_\phi^{x,y}$ generated by the variational SDE \eqref{eq:postz} is absolute continuous with respect to the path measure $\mathbb{P}^y_\theta$ induced by \eqref{eq:prior-sde}. Using the standard variational technique to apply Jensen's inequality, 
\begin{equation}
    \label{eq:jensen}
\begin{aligned}
\log p_\theta(x_{0:N}\mid y_{0:N})
&=
\log \int_\Omega
p_\theta(x_{0:N}\mid y_{0:N}, z)
\frac{d \mathbb P_\theta^y}{d \mathbb Q_\phi^{x,y}}(z)\,
\mathbb Q_\phi^{x,y}(dz) \\
&\ge \mathbb E^{\mathbb Q_\phi^{x,y}} \left[
\log p_\theta(x_{0:N}\mid z) + \log \frac{d \mathbb P_\theta^y}{d \mathbb Q_\phi^{x,y}}(z)
\right]
\\
&=
\mathbb E^{\mathbb Q_\phi^{x,y}} \left[
\log p_\theta(x_{0:N}\mid z)
\right]
-
\KL \left(\mathbb Q_\phi^{x,y}\,\|\,\mathbb P_\theta^y\right).
\end{aligned}
\end{equation}
where the last equality holds since $x_{0:N}$ and $y_{0:N}$
are conditionally independent given $z := (z_{t_k})_{k=1}^N$, so $p_\theta(x_{0:N}\mid y_{0:N}, z) = p_\theta(x_{0:N}\mid z)$. 
{The KL divergence term in the last equality of \eqref{eq:jensen} is computed analytically by \eqref{eq:variation-rn}.
The stochastic integral $\int_0^T {\alpha}_{\phi,\theta}(z_t,t) dW_t$ is a martingale with zero expectation under $\mathbb P_\theta^y$.
}

By Assumption~\ref{asmp:var_prior}, we can further simplify \eqref{eq:jensen} as \eqref{eq: condelbo}. 
By definition of $\Pi_\theta^{x,y}$, it follows that
\[
\begin{aligned}
\KL \left(\mathbb Q_\phi^{x,y}\,\|\,\Pi_\theta^{x,y}\right)
&=
\int_\Omega
\log \left(
\frac{d \mathbb Q_\phi^{x,y}}{d \mathbb P_\theta^y}(z)
\frac{d \mathbb P_\theta^y}{d\Pi_\theta^{x,y}}(z)
\right)
\,Q_\phi^{x,y}(dz) \\
&=
\mathbb E^{\mathbb Q_\phi^{x,y}} \left[
\log \frac{d \mathbb Q_\phi^{x,y}}{d\mathbb P_\theta^y}(z)
+
\log \frac{p_\theta(x_{0:N}\mid y_{0:N})}
          {p_\theta(x_{0:N}\mid z)}
\right].
\end{aligned}
\]
Rearranging the terms above, it yields that for every $(\theta,\phi)$ and training pair $(x_{0:N},y_{0:N})$,
\begin{equation}
    \label{eq: elbogap}
\begin{aligned}
\log p_\theta(x_{0:N}\mid y_{0:N})
&=
\mathcal L(\theta,\phi;x_{0:N},y_{0:N})
+
\KL \left(
Q_\phi^{x,y}\,\big\|\,\Pi_\theta^{x,y}
\right).
\end{aligned}
\end{equation}
For fixed $\theta$ and fixed data, the quantity
$\log p_\theta(x_{0:N},y_{0:N})$ is independent of $\phi$, and thus \eqref{eq:phi_opt_is_KLproj}. 
\end{proof}

We can now connect the variational approximation scheme to controlled process and $\Pi^y$ in Section 3.
Taking expectation of \eqref{eq: elbogap} over the true posterior $X \sim \Pi^y$ yields the population-level ELBO gap, $\mathcal E_y(\theta,\phi) := \mathbb E^{\Pi^y}\!\left[\log p^{\mathrm{true}}(X\mid y)\right]
-\ \mathbb E^{\Pi^y}\!\left[\mathcal L(\theta,\phi;X,y)\right]$.
With sufficient model expressivity in encoder, decoder, and the parametrization of SDE coefficients, as well as successful optimization and negligible statistical error, it is reasonable to treat this population-level ELBO gap as small. Note that we do not quantify the exact training error as universal approximation type of results, these have been discussed for some special cases in \cite{kwossek2025universal}, and classical universal approximation theorem such as \cite{hornik1989multilayer} can be adopted for the parametrization of SDE coefficients. In the following corollary, we show that the approximation error of path observables
for example a hitting event, a running cost, or a terminal-time
quantity, is controlled by the ELBO gap.  

\begin{cor}
\label{cor:gap}
Let $\Psi:\mathcal C([0,T];\mathbb R^n) \to\mathbb R$
be any bounded measurable path observable.
The true and learned posterior expectation of $\Psi$ is
\[
        \Pi^y(\Psi)
        :=
        \int_{\mathcal X_T}\Psi(x_{[0,T]})\,\Pi^y(dx_{[0,T]}), \quad \widehat\Pi_{\hat\theta}^y(\Psi)
        :=
        \int_{\mathcal X_T}\Psi(x_{[0,T]})\,
        \widehat\Pi_{\hat\theta}^y(dx_{[0,T]}),
\]
where $\widehat\Pi_{\hat\theta}^y
        =
        (\mathcal D_{\hat\theta})_\sharp\mathbb P_{\hat\theta}^y$. Assume $\Pi^y \ll \widehat\Pi_\theta^y$ on
$\mathcal C([0,T];\mathbb R^n)$.
 If the learned parameters $(\hat\theta,\hat\phi)$ achieve conditional population ELBO gap at most $\varepsilon_{\mathrm{opt}}$, i.e.
$\mathcal{E}_y(\hat\theta,\hat\phi)\le \varepsilon_{\mathrm{opt}}$, then 
\[
\left|
        \Pi^y(\Psi)
        -
        \widehat\Pi_{\hat\theta}^y(\Psi)
\right|
\leq \|\Psi\|_\infty\sqrt{2\varepsilon_{\mathrm{opt}}}.
\]
\end{cor}

\begin{proof}
Substitute
$\mathbb E^{\Pi^y}[\log p^{\mathrm{true}}(X\mid y)-\log p_\theta(X\mid y)]
=\mathrm{KL}(\Pi^y\|\widehat\Pi_\theta^y)$ into the definition of $\mathcal E_y(\theta,\phi)$, we get 
\begin{equation}
\label{eq:gap-two-kl}
\mathcal E_y(\hat\theta, \hat\phi)
=\mathrm{KL}\left(\Pi^y\,\|\,\widehat\Pi_{\hat\theta}^y\right)
+\mathbb E^{\Pi^y} \left[\mathrm{KL}(\mathbb Q_{\hat\phi}^{X,y}\,\|\,\Pi_{\hat\theta}^{X,y})\right].
\end{equation}
$\mathcal E_y$ vanishes when $\widehat\Pi_\theta^y=\Pi^y$ and
$\mathbb Q_\phi^{X,y}=\Pi_\theta^{X,y}$ $\Pi^y$-a.s.
Both terms in \eqref{eq:gap-two-kl} are nonnegative and thus are bounded by $\varepsilon_{\mathrm{opt}}$ individually. In particular, the first term in the last equality shows the approximation error of the decoded generative model from the true posterior distribution; the second term shows the expected variational gap between $\mathbb{Q}^\cdot$ and $\mathbb{P}^\cdot$ during training.

By Pinsker's inequality, for a bounded measurable $\Psi$,
    \[
\left|
        \Pi^y(\Psi)
        -
        \widehat\Pi_{\hat\theta}^y(\Psi)
\right|
\leq
        2\|\Psi\|_\infty
        \,
        \|\Pi^y-\widehat\Pi_{\hat\theta}^y\|_{\mathrm{TV}}
\leq
        \|\Psi\|_\infty
        \sqrt{
        2\KL\!\left(\Pi^y\,\middle\|\,\widehat\Pi_{\hat\theta}^y\right)
        } \leq \|\Psi\|_\infty\sqrt{2\varepsilon_{\mathrm{opt}}}.
\]
\end{proof}

Suppose $\widehat\Pi_{\hat\theta}^y$ is the law
of an It\^o diffusion sharing the signal's diffusion coefficient $\sigma$, with drift
$b_{\hat\theta}(t,X_t;y_{[0,t]})$. By Theorem~\ref{thm: pathcontrol}, $\Pi^y$ is the same diffusion
driven by the \emph{optimal} feedback drift
$b^\star(t,x;y_{[0,t]}):=\beta(t,x)+\sigma(t,x)\,\alpha^\star(t,x;y_{[0,t]})$.
By Girsanov's theorem,
\begin{equation}
\label{eq:girsanov-main-drift}
\KL\left(\Pi^y\,\|\,\widehat\Pi_{\hat\theta}^y\right)
=\KL\left(\nu_0^y\,\|\,\widehat\nu_{\hat\theta,0}^y\right)
+\tfrac12\,\E^{\Pi^y}\!\left[\int_0^T
\big\|\sigma(t,X_t)^{-1}\left(b^\star-b_{\hat\theta}\right)(t,X_t;y_{[0,t]})\big\|^2\,dt\right]
\;\le\;\varepsilon_{\mathrm{opt}},
\end{equation}
where $\nu_0^y,\widehat\nu_{\hat\theta,0}^y$ are the initial marginals of $\Pi^y$ and
$\widehat\Pi_{\hat\theta}^y$. Thus, the pathwise ELBO drives the generator drift to the optimal
control of Theorem~\ref{thm: pathcontrol}, and trains the sampler to approximate a controlled diffusion that induces $\Pi^y$.

\subsection{Inference: subproblems of posterior law}
\label{sec:sampler}

Recall from Section~\ref{sec: theory} that by letting the terminal time range over the time grid $0=t_0<t_1<\dots<t_N=T$, there are controlled processes  producing a sequence of subproblems $\Pi^{y,t_k}_{0}:=\Law\!\left(X_{[0,t_k]}\mid Y_{[0,t_k]}=y_{[0,t_k]}\right)$, for $k=0, \ldots, N$. 
Each $\Pi^{y,t_k}_{0}$ is realized by a controlled diffusion on $[0,t_k]$ whose feedback drift is the score $-\sigma\sigma^\top\nabla_x\widetilde V^y_{t_k}$ by Theorem~\ref{thm: pathcontrol}.
However, when simulating the controlled process by Euler-Maruyama, it integrates the learned generator drift over
$[t_k,t_{k+1}]$, producing the predictive law $\Law(Z_{[0,t_{k+1}]}\mid y_{0:k})$. This is the predict step. To achieve the correct filtering distribution, we need an update
step that applies the incremental log-likelihood at the
newly arrived observation $y_{k+1}$, turning the predictive law into the filtering law
$\Law(Z_{[0,t_{k+1}]}\mid y_{0:k+1})$. 
In particular, with step $\Delta t$ and $\Delta W_k\sim\mathcal N(0,\Delta t\,I)$,
\begin{align}
z_{k+1}^-
&=
z_k+\mu_\theta\!\left(t_k,z_k;y_{0:k}\right)\,\Delta t
+g_\theta(t_k,z_k)\,\Delta W_k,
\label{eq:em-predict}\\[2pt]
z_{k+1}
&=
z_{k+1}^-
+a_\theta(t_{k+1},z_{k+1}^-)\,
   \nabla_{z}\log p_\theta\!\left(y_{k+1}\mid z_{k+1}^-\right)\,\Delta t,
\qquad a_\theta:=g_\theta g_\theta^\top.
\label{eq:em-update}
\end{align}
After decoding, this reproduces $\Pi^{y,t_k}_{0}$ at \emph{every} $t_k$. More specifically, since
\[
x_{t_i}=\mathcal D_\theta(z_{t_i}),
\qquad
y_{t_i}=h(x_{t_i})+\epsilon_{t_i},
\qquad
\epsilon_{t_i}\sim\mathcal N(0,R),
\]
then, up to an additive constant,
\begin{equation}
\label{eq: obsll}
\log p_\theta\!\left(y_{t_i}\mid z_{t_i}\right)
=
-\tfrac12\left(y_{t_i}-h(\mathcal D_\theta(z_{t_i}))\right)^\top R^{-1}
       \left(y_{t_i}-h(\mathcal D_\theta(z_{t_i}))\right).
\end{equation}

\begin{remark}[Comparison with particle-based methods]
In our method, no knowledge of the signal dynamics is needed during training. That is, the model does not use the prior measure $\mu$ induced by \eqref{eq: state}, but learns a map from noisy observation paths to the posterior path law. Classical particle filters or SMC, by
contrast, require the transition kernel to propagate and weight particles, and
their accuracy is sensitive to misspecification of the prior dynamics. 

The two formalisms also treat
the observation likelihood update differently.
SMC approximates the marginal distributions at each time and update the distribution by multiplicative importance weights. This often requires resampling to counter weight
degeneracy phenomenon. In our method, the running-cost potential $\widehat G(t,\cdot;y)$ of the pathwise Zakai
equation (Proposition~\ref{prop:path-zakai}, Theorem~\ref{prop: control}) is the
continuous-time analog of the incremental log-likelihood $\log p(y_t\mid Z_t)$. Once trained, the model produces i.i.d. 
posterior path samples by a single forward simulation of \eqref{eq:em-predict}-\eqref{eq:em-update} with the decoder $\mathcal D_\theta$.
\end{remark}

\section{Experiments}
\label{sec: experiments}

In this section, we present several examples on the inference of the stochastic dynamics from noisy, non-linear, and possibly missing  observations, using our method. The implementation is done in Python. \footnote{  \url{https://github.com/nicoletyang/FilteringSDE.git}}

\subsection{Double-well equation: Computation of path functionals}
We consider a double-well equation with stochastic parts. For $t \in [0,T]$,
\begin{equation}
\label{eq:dw-1}
    dX_t = -4 X_t (X_t^2 - 1) dt + \beta dW_t, \quad X_0 \sim N(0, I).
\end{equation}
The deterministic part of the system has two stable fixed points (wells) $X_t = 1$ and $X_t = -1$, and an unstable fixed point $X_t = 0$. We set $\beta =1$ so that there are transitions between the stable states. The stationary density is bimodal, $p_X \propto e^{-(y^4 - 2y^2)/\beta^2}$.

\begin{figure}
    \centering
\includegraphics[width=1\linewidth]{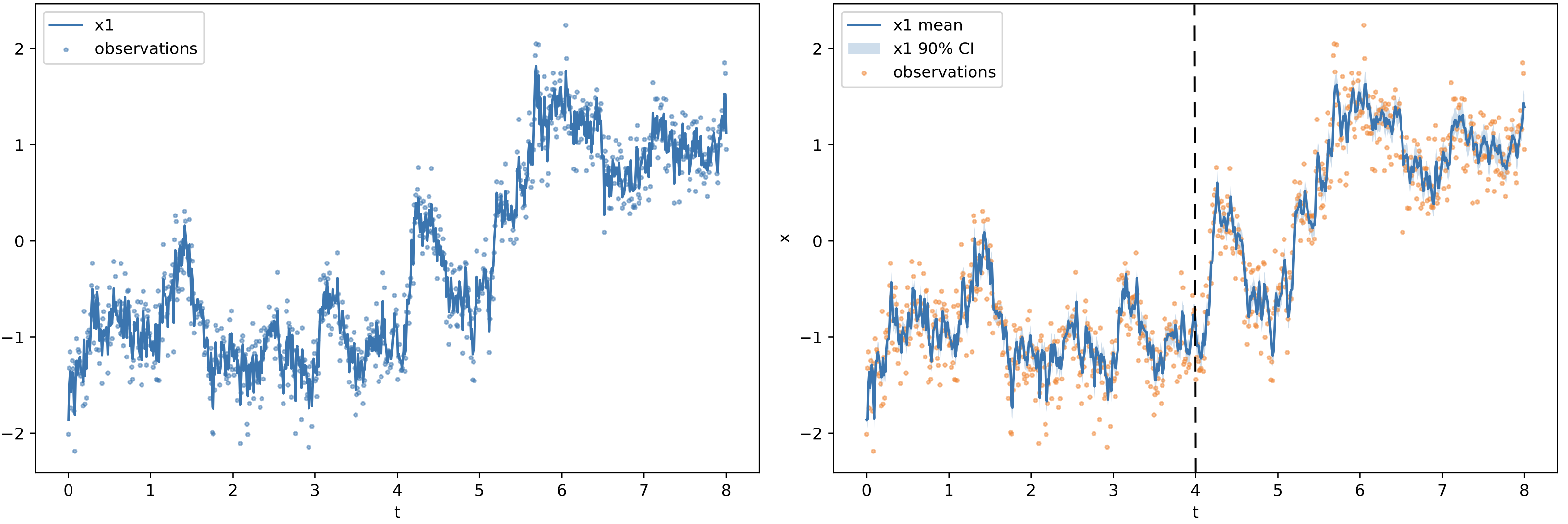}
    \caption{\textbf{Stochastic double well equation}. The model is trained on synthetic data from \eqref{eq:dw-1} on time horizon $[0,4]$ and perform inference on test noisy observation paths on time horizon $[0,8]$. Left: True test trajectories and noisy observations; Right: Estimated trajectories and the 90\% confidence interval (CI) with only the noisy observations available during inference.}
    \label{fig:placeholder}
\end{figure}

We simulate a batch of training and test trajectories synthesized from \eqref{eq:dw-1} and generate noisy observations $y_t = x_t + \varepsilon_t$ with $\varepsilon_t \sim \mathcal{N}(0,\sigma^2 I)$. We train the model following Section~\ref{sec: nsde} and the optimization of ELBO is done with Adam and a decaying learning rate. At evaluation time, we draw many posterior sample paths by running \eqref{eq:em-update}.
Then, we use these samples to 
assess path estimation quality such as RMSE, Wasserstein-1 distance; demonstrate marginal density via histograms of the estimated paths and the true simulated paths; 
and more interestingly, estimate mean dwell times, a functional $\tau(\cdot)$ of the path posterior law. 

\begin{figure}[htp] 
\centering{
\includegraphics[scale=0.4]{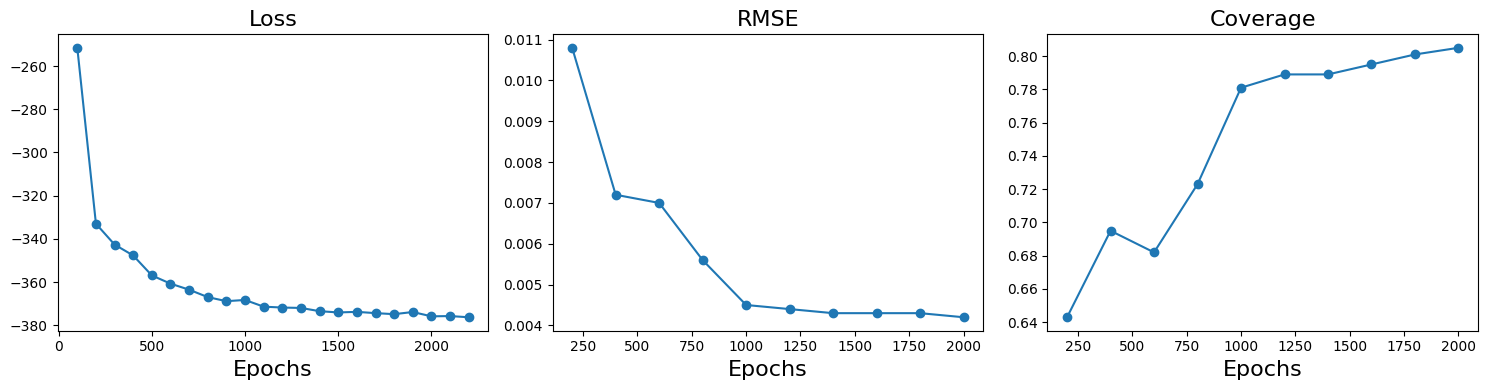}
}
\caption{Estimated trajectories from stochastic double well equation. Left to right: Loss every 100 epochs, Dwell time RMSE and 90\% coverage every 200 epochs. The inferred trajectories from our method can capture long-time, metastable behavior of the underlying dynamics, verifying the efficient learning of the posterior path measure. 
}
\label{fig:doublewell}
\end{figure}

In particular, we observe a one-dimensional state path $x_1,\dots,x_{N}$ at times
$0 = t_1 < \dots < t_{N} = T$. Fix a threshold $c \in \mathbb{R}$ (in our case $c = 0$) and 
we consider the dwell-time functional
$\tau(x) \approx \frac{1}{N}\sum_{k=1}^{N} \mathbf{1}\{x_{t_k} < 0\}$.
We estimate the dwell times of the test data (ground truth paths) and our estimated paths by Monte-Carlo simulation and present the RMSE in Figure~\ref{fig:doublewell}.

For each Monte Carlo sample
$b$, let $q^{(b)}_{0.05}$ and $q^{(b)}_{0.95}$ be the $5\%$ and $95\%$
quantiles of $\tau(x^{(\ell)}_{b})$. 
The 90\%-coverage
metric for occupation time is then
$\frac{1}{B}\sum_{b=1}^B 
\mathbf{1}\Big\{\tau^{\text{true}}_b \in [q^{(b)}_{0.05},\,q^{(b)}_{0.95}]\Big\},$
which estimates the empirical coverage probability of $90\%$ posterior
intervals for the expected occupation time across trajectories.

\subsection{Lorenz-63: Low data requirement and comparison with particle-based methods}

\begin{figure}[t]
  \centering
\includegraphics[trim={0 0.8cm 0 0},clip, width=1\linewidth]{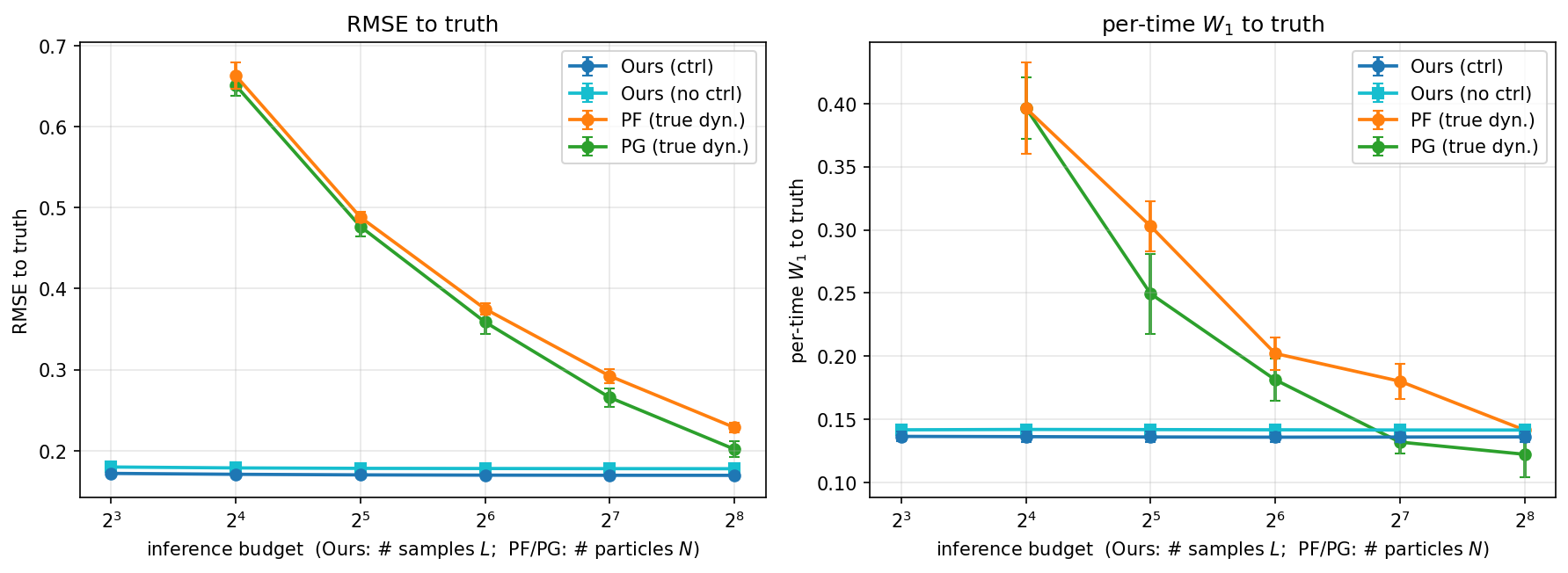}
  \caption{
  \textbf{Stochastic Lorenz-63.}
  Comparison of posterior accuracy as a function of the inference budget ($\log_2$ axis), averaged over five seeds; error bars show one standard error. The budget is represented by the number of posterior path samples $L$ for our method and the number of particles $N$ for the particle filter (PF) and particle Gibbs (PG); PF and PG are additionally given the true dynamics~\eqref{eq: lorenz63}, while our method learns them from data. 
We use the observation model $Y_t = \arctan(X_t) + \varepsilon_t$, $\varepsilon_t \sim \mathcal{N}(0,\sigma_{\mathrm{obs}}^2 I_3),$ with PF/PG applying likelihood updates only at observed indices while propagating true dynamics \eqref{eq: lorenz63}.   Left: RMSE of the posterior mean to the ground-truth path. Right: per-time Wasserstein-1 distance to the ground-truth path, averaged over time, batch, and coordinate.
  }
  \label{fig:pf_pg_vs_latentsde}
\end{figure}

The Lorenz attractor exhibits chaotic behavior on the strange attractor. With the additional stochasticity, we take $\theta = (10, 28, 8/3)$ and $\beta = (.1, .28, .3)$, 
\begin{equation}
\label{eq: lorenz63}
    \begin{aligned}
    dX_1(t) &= \theta_1(X_2(t) - X_1(t)) dt + \beta_1 X_1(t) dW_t, \\
    dX_2(t) &= (\theta_2 X_1(t) - X_1(t) X_3(t)) dt + \beta_2 X_2(t) dW_t, \\
    dX_3(t) &= (X_1(t) X_2(t) - \theta_3 X_3(t)) dt + \beta_3 X_3(t) dW_t.
\end{aligned}
\end{equation}
The stationary distribution exists, see \cite{keller1996attractors} for analysis on stochastic Lorenz system. 
We use this as a testbed to compare with traditional particle-based methods, Particle Filters (PF) and Particle Gibbs (PG), and present the inferred trajectories in Figure~\ref{fig:pf_pg_vs_latentsde}. PG targets the posterior distribution by repeated conditional particle filter updates and is asymptotically exact. In sparse observation regimes, conditional SMC may suffer from weight and path degeneracy and the PG kernel may mix slowly: trajectories remain stuck near the reference path unless the number of particles is taken large enough. Ancestor sampling helps but does not remove the fundamental issue that the posterior is far from the prior in sparse chaotic settings. Note that particle-based methods require the signal dynamics. In the below comparison results, the correct Lorenz system is given in particle-based methods, while in our method the posterior dynamics is learned directly from the data, without knowledge of the prior dynamics. 

{Our method attains a reasonable accuracy (mean RMSE $= 0.1719$ and mean $W_1 = 0.1364$ over 5 seeds) with 8 samples and the accuracy increases slightly as the number of samples grows. With 256 samples, we get mean RMSE $= 0.169550$ and mean $W_1 = 0.136047$ over 5 seeds. The particle baselines require one to two orders of magnitude more particles to approach it. On $W_1$ our method dominates PF throughout and is overtaken by PG when the number of particles is above 128, reflecting the tighter trajectory tracking that the exact dynamics afford PG. In our method, we compare two parametrizations, the full $w_\theta(t,z; y_{[0,t]}) + g_\theta u_\phi(t,z; y_{[0,t]})$ as the drift, versus a single network $\mu_\theta(t,z; y_{[0,t]})$ as the drift. The structure $w_\theta + g_\theta u_\phi$ gives marginally better result and may serve as a useful inductive bias to improve the expressivity of the learned drift.}

\subsection{Lorenz-96: Robustness towards sparse \& missing observations}

For a higher-dimensional example, we perform the estimation tasks on noisy observations from a stochastic 15-dimension Lorenz-96 \cite{lorenz1996predictability} system, where for $i = 1, \ldots, 15$, $F = 8$,
\[
dX_{t,i} = \left( (X_{t, i+1} - X_{t, i-2}) X_{t, i-1} - X_{t,i} + F \right) dt + b_i X_{t, i} dW_{t,i}.
\]
Recall that the dynamics is not used during training or inference, it is only used to simulate the data. Here we simulate training data from [0,2], with 200 time steps.

\begin{figure}[htp] 
\centering{
\includegraphics[scale=0.4]{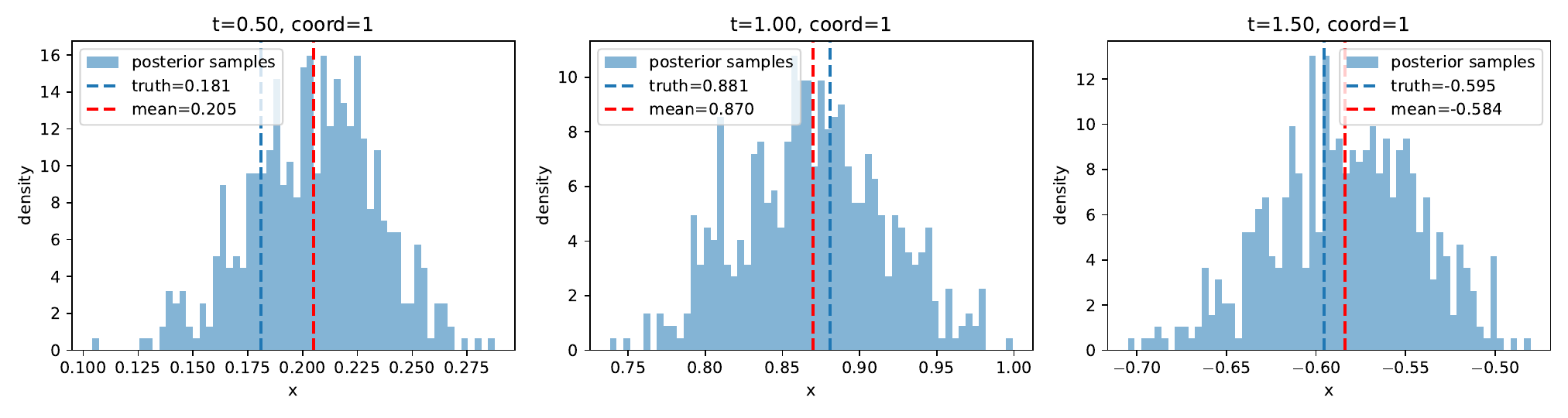}
\includegraphics[scale=0.32]{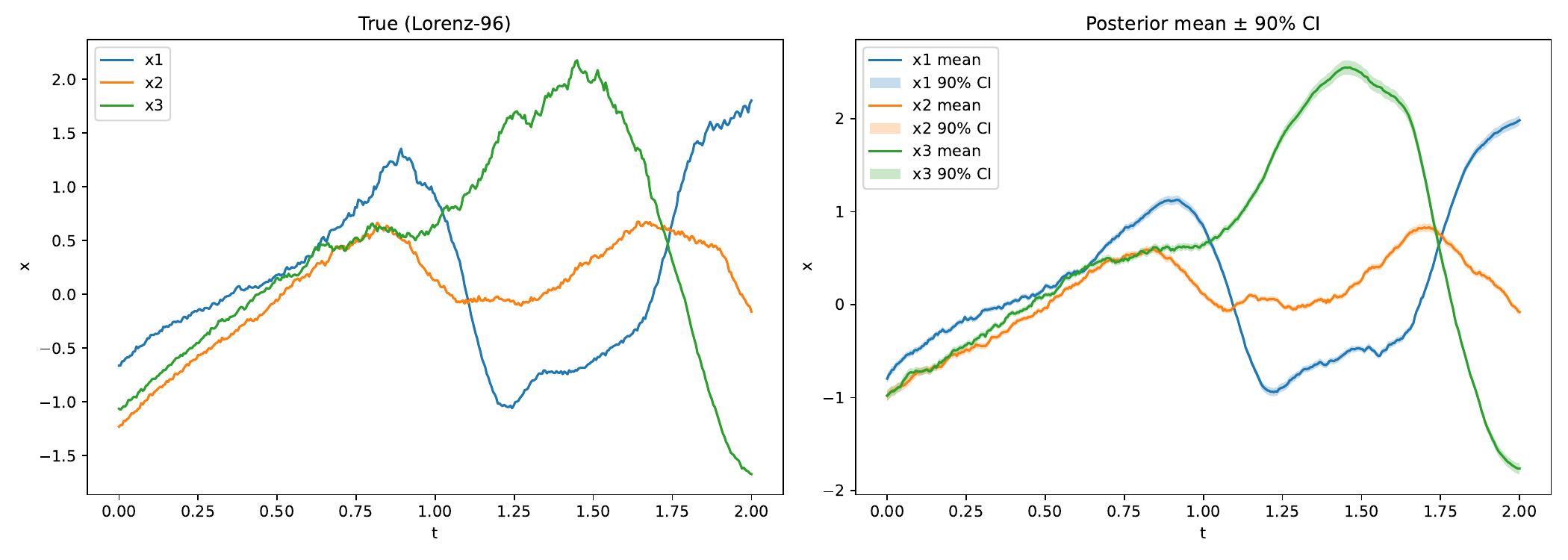}
}
\caption{\textbf{Stochastic Lorenz-96 equation.} Estimated trajectories from $15$-dimensional stochastic Lorenz-96 equation. The model is trained for time $[0,2]$ with only observation from test dataset available. The observation model is $y_t = Tanh(x_t)+ N(0, \sigma^2), \sigma=0.15$. Top: Comparison of mariginal distributions at time 0.5, 1, and 1.5 for the first dimension. Bottom: True (left) and Inferenced (right) trajectories of the first 3 dimensions. 90\% confidence intervals are also presented for the inferenced trajectories.}
\label{fig:lorenz96lin}
\end{figure}

In Figure~\ref{fig:lorenz96lin}, we show the inference result in terms of the marginal distribution at a few time points as well as the inferred curves with the confidence interval.

\begin{figure}[htp] 
\centering{
\includegraphics[scale=0.3]{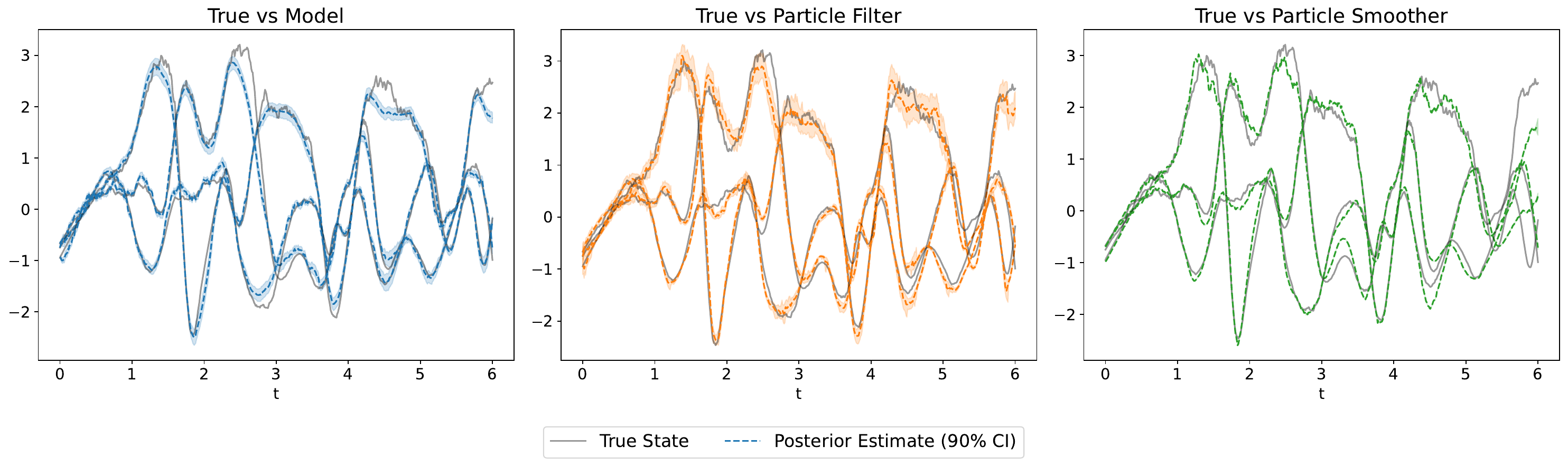}
}
\caption{Estimated trajectories of the first 3 dimensions from $15$-dimensional stochastic Lorenz-96 equation. The model is trained for time $[0,3]$ with 20\% observation randomly masked. During inference, only the noisy observation from test dataset available, (random) 20\% of the observation time is missing. The observation model is $y_t = arctan(x_t)+ N(0, \sigma^2), \sigma=0.15$. 90\% confidence intervals are also presented for the inferenced trajectories.}
\label{fig:lorenz96nonlin}
\end{figure}

In Figure~\ref{fig:lorenz96nonlin}, we test the effect of missing observations, comparing against particle-based inference baselines under the same observation model and missing rate. In particular, to model incomplete observations, we introduce a binary mask $m_t = (m_{t,1},\dots,m_{t,D}) \in \{0,1\}^D$,
where $m_{t,d} = 1$ if coordinate $d$ is observed at time $t$, and $m_{t,d} = 0$ if not.
In experiments, the mask is sampled independently for each time and coordinate $
    m_{t,d} \sim \mathrm{Bernoulli}(1-r)$,
where $r \in [0,1]$ is the prescribed missing rate.
Equivalently, the observation input to the model is $\bar{y}_t = m_t \odot y_t$,
where $\odot$ denotes elementwise multiplication. 
The observation encoder receives both the masked observations and the mask itself to distinguish 
the true observed value being close to zero versus the masked/missing observation that is set to zero. In particular, the observation-context encoder here is a masked causal self-attention layer
$E_n = \Psi_{\theta}^{\mathrm{attn}}
    \bigl(\widetilde Y_0,M_0,\ldots,\widetilde Y_n,M_n\bigr)$. 
The approximate posterior over the initial latent state uses both an $x$-context encoder and an observation-context encoder, and the prior drift depends only on the observation-context encoder.

We further present in Table~\ref{tab:missing-rate-metrics} the evaluation metrics on the path inference error across different missing rates and compare with particle-based baseline methods. As the missing rate grows from $10\%$ to
$50\%$, our path-inference performance degrades slightly, whereas the particle baselines
degrade markedly faster. 

\begin{table}[t]
\centering
\caption{\textbf{$15$-dimensional stochastic Lorenz-96 equation.} Evaluation metrics under different missing rates over 10 seeds. The model is trained for time $[0,3]$ with 20\% observation randomly masked. During inference, only the noisy observation from test dataset available, with 10\%-50\% of the observation time missing. PF and PG are implemented with 512 particles. Sample size of our approach is 64. RMSE and $W_1$ are recorded as mean $\pm$ standard error, lower is better.}
\label{tab:missing-rate-metrics}
{\scriptsize
\begin{tabular}{@{}c c c c@{}}
\toprule
\textbf{Missing Rate} & \textbf{Method} & \textbf{RMSE} & $\mathbf{W_1}$ \\
\midrule
{0.1}
& Ours & $0.141{\scriptstyle\pm0.002}$ & $0.104{\scriptstyle\pm0.002}$ \\
& PF   & $0.223{\scriptstyle\pm0.004}$ & $0.128{\scriptstyle\pm0.010}$ \\
& PG   & $0.219{\scriptstyle\pm0.005}$ & $0.106{\scriptstyle\pm0.006}$ \\
\midrule
{0.2}
& Ours & $0.148{\scriptstyle\pm0.001}$ & $0.109{\scriptstyle\pm0.003}$ \\
& PF   & $0.237{\scriptstyle\pm0.004}$ & $0.178{\scriptstyle\pm0.028}$ \\
& PG   & $0.237{\scriptstyle\pm0.004}$ & $0.136{\scriptstyle\pm0.015}$ \\
\midrule
{0.3}
& Ours & $0.159{\scriptstyle\pm0.002}$ & $0.117{\scriptstyle\pm0.002}$ \\
& PF   & $0.281{\scriptstyle\pm0.004}$ & $0.182{\scriptstyle\pm0.016}$ \\
& PG   & $0.280{\scriptstyle\pm0.004}$ & $0.153{\scriptstyle\pm0.020}$ \\
\midrule
{0.4}
& Ours & $0.175{\scriptstyle\pm0.002}$ & $0.124{\scriptstyle\pm0.002}$ \\
& PF   & $0.326{\scriptstyle\pm0.003}$ & $0.177{\scriptstyle\pm0.017}$ \\
& PG   & $0.322{\scriptstyle\pm0.006}$ & $0.175{\scriptstyle\pm0.023}$ \\
\midrule
{0.5}
& Ours & $0.195{\scriptstyle\pm0.002}$ & $0.137{\scriptstyle\pm0.004}$ \\
& PF   & $0.377{\scriptstyle\pm0.004}$ & $0.237{\scriptstyle\pm0.021}$ \\
& PG   & $0.378{\scriptstyle\pm0.003}$ & $0.242{\scriptstyle\pm0.026}$ \\
\bottomrule
\end{tabular}
}
\end{table}

\subsection{Real dataset: MuJoCo Hopper}

We use the Hopper physics simulation from the DeepMind Control Suite \cite{tassa2018deepmind}. Each trajectory is generated by sampling a random initial configuration and velocity, then rolling out the deterministic simulator under a zero-action policy. Let $x(t)\in\mathbb{R}^{14}$ denote the simulator state at time $t$, consisting of concatenated generalized positions and velocities. We generate $10{,}000$ trajectories of length $T=100$ on an evenly spaced grid $0=t_0<\cdots<t_{T-1}=1$ and standardize each dimension using training-set statistics. Observations are corrupted by additive Gaussian noise $y_{t_i} = x_{t_i} + \varepsilon_{t_i}, \, 
\varepsilon_{t_i}\sim \mathcal{N}(0,\sigma_{\mathrm{obs}}^2 I_{14})$. We use a GRU encoder for the noisy observations here.

As an ablation study with training marginal distributions,
we compare against a probabilistic GRU autoregressive model that defines a
causal factorization as
\begin{equation}
x_i \mid x_{i-1}, y_{\le i} 
\sim
\mathcal{N} \left(\mu_\theta(h_i),\mathrm{diag}(\sigma_\theta^2(h_i))\right), \quad p_\theta(x_{0:N}\mid y_{0:N})
=
p_\theta(x_0\mid y_0)
\prod_{i=1}^{N}
p_\theta(x_i \mid x_{i-1}, y_{\le i}),
\label{eq:gru_ar_factorization}
\end{equation}
where the hidden state $h_i$ is updated by a GRU using only information
available up to time $t_i$.
The output heads are
$\mu_\theta(h_i)=W_\mu h_i+b_\mu$, and $\sigma_\theta(h_i)=\mathrm{softplus}(W_\sigma h_i+b_\sigma)+\sigma_{\min}$.
We train by maximizing the conditional log-likelihood, $\max_\theta \sum_{i=0}^{N}
\log p_\theta(x_i \mid x_{i-1}, y_{\le i}),$
and at test time the model conditions only on noisy observations $y$ and masks, while $x$ is used only for evaluation.

\begin{figure}
    \centering
\includegraphics[width=0.8\linewidth]{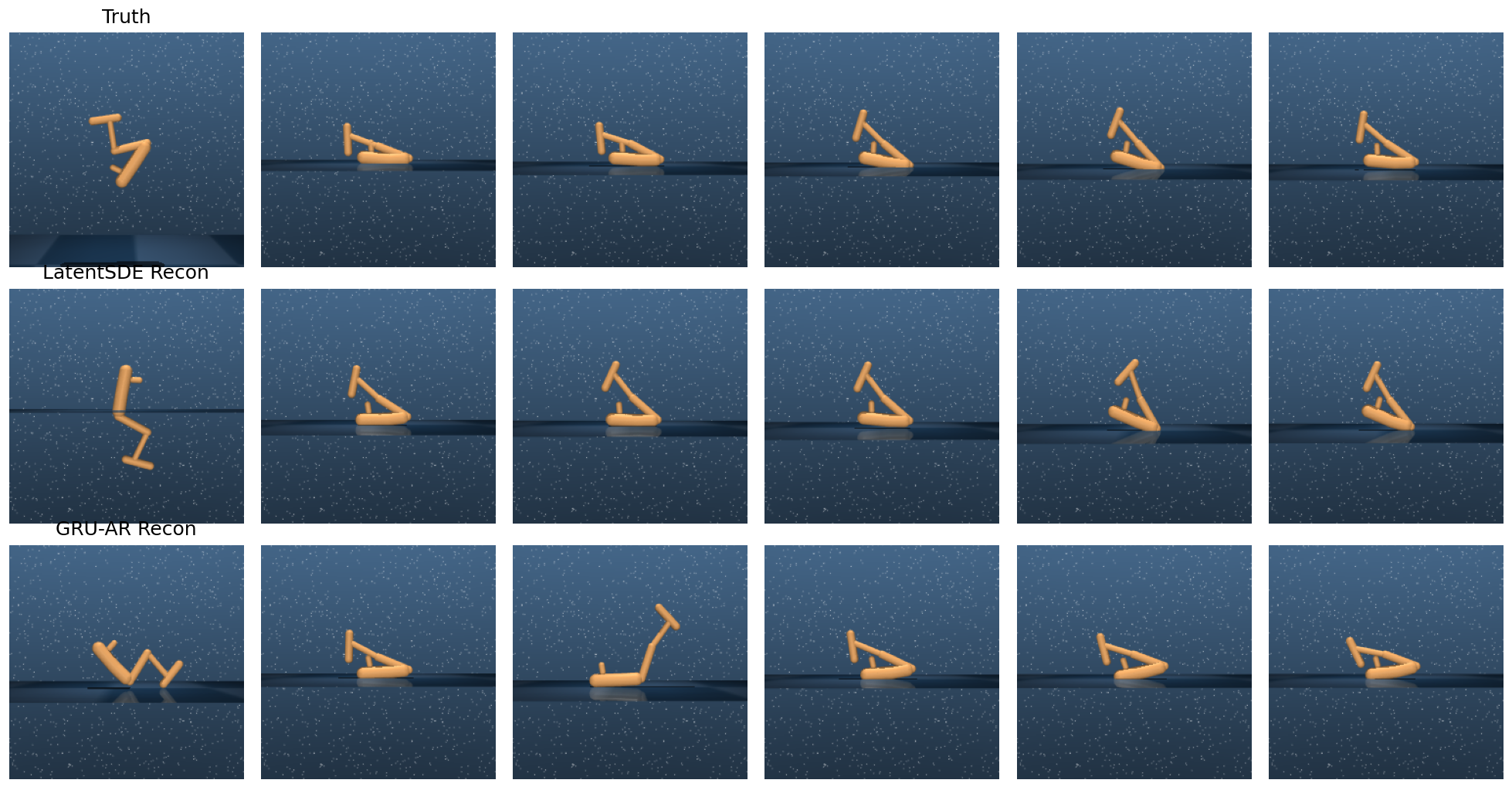}
\includegraphics[width=0.8\linewidth]{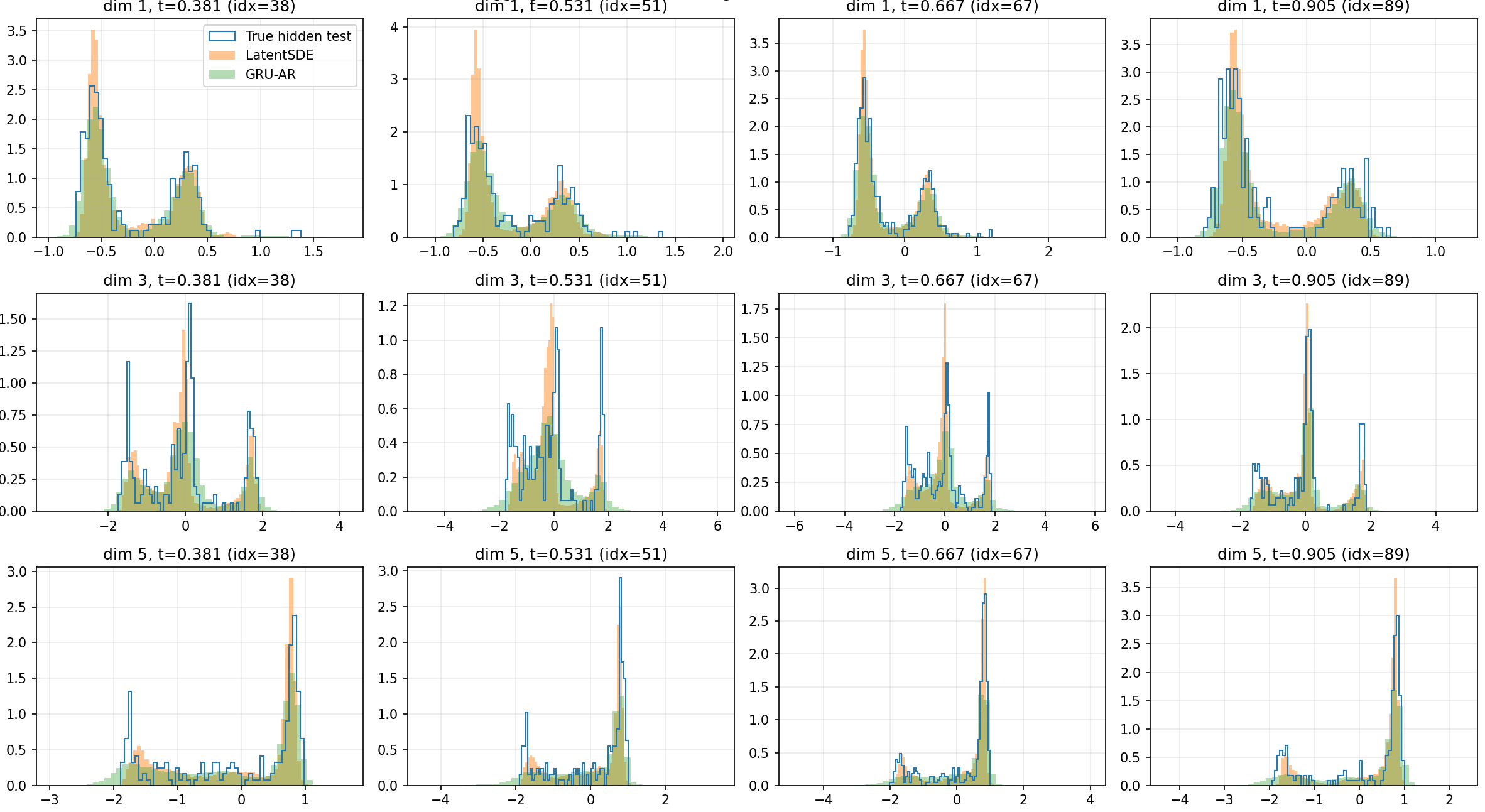}
    \caption{\textbf{Hopper physics simulation}. We compare the assimilated results from our conditional latent SDE approach against the GRU autoregressive model. Left: Inferred trajectories. Top row represents the true trajectory of the Hopper in its physical environment. Second and third row show the inferred frames from our conditional latent SDE method vs GRU-AR, respectively. Right: We take dimension 1,3,5 which has an obvious multimodal pattern in the distribution and plot the corresponding histogram at several time slices (0.38, 0.53, 0.67, 0.9). We report the pathwise Wasserstein-2 ($W_2$) distance between the inferred trajectories and the true test trajectories over five random seeds as mean $\pm$ standard error. For conditional latent SDE, the $W_2$ evaluated on the entire time hoirzon and the $W_2$ evaluated on the missing time windows are 
 $29.13 \pm 3.04$ and $11.51 \pm 3.05$, respectively; for GRU-AR, the $W_2$ evaluated on the entire time hoirzon, and the $W_2$ evaluated on the missing time windows are $30.33 \pm 3.13$ and $14.24 \pm 3.28$, respectively.}
    \label{fig:hopper}
\end{figure}

Empirically, GRU-AR can be highly competitive on short-horizon denoising when the mapping $y_{\le t}\mapsto x_t$ is smooth and unimodal; however, it estimates marginal state locally in time but does not define a continuous-time latent process and does not model posterior path dependence beyond the recurrent summary $h_i$. We see that in Figure~\ref{fig:hopper}, with multi-modal and irregularly sampled data, the conditional latent SDE encourages temporal coherence and can better capture the different modes in the data. 

\section{Discussion and Future Work}

In this paper, we derive an optimal control formulation of pathwise filtering and develop an observation-conditioned neural SDE that induces the posterior path law. This provides an amortized approximation to the Bayesian posterior on trajectories.
The conditional latent SDE framework offers several advantages over traditional particle-based methods, 
when (i) the transition mechanism (drift and diffusion) is unknown or only partially specified, (ii) repeated assimilation across many sequences is required, and (iii) uncertainty quantification and computation of functionals over trajectories are required under possibly sparse observation. Furthermore, the continuous‑time formulation helps the algorithm to deal with randomly missing observations during inference.
In data assimilation language, the method here achieves both parameter and state estimation. Rather than approximating the marginal posterior distribution with particle methods, here we approximate the prior by the generative model and correct the sample path by the optimal control solution developed upon Zakai equation \cite{zakai1969optimal}, alternatively, one can understand our method as learning the `posterior' SDE that induce the desired posterior path measure. This perspective opens up several promising directions for future work.

First, the continuous-time formulations are naturally suited to missing or irregularly sampled data, but theoretical guarantee and modeling choices remain important. In the online setting, one may seek explicit stability or robustness statements (\cite{crisan2010approximate}) for the learned conditional kernel $y\mapsto \PP_\theta^y$.
Second, the resulting system can be viewed as a McKean-Vlasov type, and one can study training and approximation via propagation-of-chaos or mean-field limits, with techniques in \cite{arous1999increasing, yang2023relative} for example. It would be valuable to quantify approximation errors in terms of statistical error from finite data and optimization error from the neural parameterization.

This problem can also be understood as an entropy regularized problem on path space relative to a reference diffusion. Thus, it shares similarity with Schr\"odinger Bridge (SB)problems, where SB uses hard marginal constraint at initial and terminal times, while the problem here uses soft constraints thought the likelihood potential conditional on the noisy observation paths. We can construct a framework for data assimilation that unifies different data assimilation methods from a generalized Schr\"odinger Bridge \cite{chen2016relation, deng2024reflected, chen2023provably} or
stochastic interpolants \cite{chen2024probabilistic, albergo2023stochastic} perspective, and consider transformer-based architecture such as in \cite{chang2025dual}. 
Computation-wise, as we have a control network and an attention-based observation encoder, there will be a quadratic cost added to the memory and time complexity, comparing to the constant memory and $O(L \log L)$ time complexity as in Latent SDE \cite{li2020scalable} without noisy observation, where $L$ the the number of steps used in a fixed-step solve. This may result in challenges with very long sequence of data. We may consider signature-based approach or rough SDEs such as in \cite{morrill2021neural}.

\section*{Acknowledgments}
I am grateful to Professor Lars Ruthotto and Professor Tomoyuki Ichiba for many stimulating discussions on data assimilation and nonlinear filtering.

\newpage
\bibliographystyle{alpha}
\bibliography{sample}

\begin{appendix}

\end{appendix}
\end{document}